\documentclass[11pt]{article}
\usepackage{mathrsfs}
\usepackage{amsmath}
\usepackage{amsfonts}
\usepackage{amssymb, amsmath, cite}
\usepackage{hyperref}
\usepackage{color}

\newtheorem{theorem}{Theorem}
\newtheorem{lemma}{Lemma}
\newtheorem{proposition}{Proposition}
\newtheorem{remark}{Remark}

\newtheorem{corollary}{Corollary}

\newcommand\be{\begin{equation}}
\newcommand\ee{\end{equation}}
\newcommand\ber{\begin{eqnarray}}\newcommand\bea{\begin{eqnarray}}
\newcommand\eer{\end{eqnarray}}\newcommand\eea{\end{eqnarray}}
\newcommand\berr{\begin{eqnarray*}}
\newcommand\eerr{\end{eqnarray*}}
\newcommand\ba{\begin{array}}
\newcommand\ea{\end{array}}

\def\Xint#1{\mathchoice
  {\XXint\displaystyle\textstyle{#1}}%
  {\XXint\textstyle\scriptstyle{#1}}%
  {\XXint\scriptstyle\scriptscriptstyle{#1}}%
  {\XXint\scriptscriptstyle\scriptscriptstyle{#1}}%
  \!\int}
\def\XXint#1#2#3{{\setbox0=\hbox{$#1{#2#3}{\int}$}
  \vcenter{\hbox{$#2#3$}}\kern-.5\wd0}}

\def\dashint{\Xint-}

\newcommand{\lm}{\lambda}\newcommand{\bfR}{\mathbb{R}}\newcommand{\pa}{\partial}

\newcommand\re{\mathrm{e}}
\newcommand\ri{\mathrm{i}}
\newcommand{\ud}{\mathrm{d}}
\newcommand{\nm}{\nonumber}\newcommand{\nn}{\nonumber}
\newcommand{\ds}{\displaystyle}

\newcommand{\itr}{\int_{\mathbb{R}^2}}
 
\newcommand{\vep}{\varepsilon}
\def\Tr{{\rm Tr\, }}

\title{Non-topological Vortex Configurations  in the   ABJM Model\footnote{This work is supported by PRIN12:  "Variational and Perturbative Aspects in Nonlinear Differential Problems" and by  FIRB project: "Analysis and Beyond".}}

\author{Xiaosen Han$^{a,b}$ \quad   \quad Gabriella Tarantello$^{a}$
   \\
\small {$^a$ Dipartimento di Matematica, Universita degli Studi di Roma ¡±Tor Vergat\`{a}¡±,}\\
\small{Via della Ricerca Scientifica, 00133 Rome, Italy}\\
\small {$^b$ Institute of Contemporary Mathematics, School of Mathematics,  Henan University}\\
\small{Kaifeng  475004, PR China}}
\date{}

\begin{document}
\maketitle
\begin{abstract}
 In this paper we study the existence of  vortex-type solutions  for a system of  self-dual equations deduced from the mass-deformed  Aharony--Bergman--Jafferis--Maldacena (ABJM) model.
 The governing equations, derived by Mohammed, Murugan, and Nastse under suitable ansatz involving fuzzy sphere matrices,  have the   new feature that
 they can support only non-topological vortex solutions. After transforming  the self-dual equations into a nonlinear elliptic
   $2\times 2$ system we  prove first an existence result by means of a perturbation argument based on a new
 and appropriate scaling for the solutions. Subsequently, we prove a more complete  existence result  by using a dynamical analysis  together with a blow-up argument.
 In this way we  establish  that,  any positive energy level is attained by a 1-parameter family of vortex solutions which also correspond to  (constraint) energy minimizers.
 In other words, we register the exceptional fact in a BPS-setting that, neither a  ``quantization'' effect nor an  energy gap is induced  upon the system by the rigid
 ``critical'' coupling of the self-dual regime.

\medskip

\medskip

 \end{abstract}
 \section{Introduction}
 \setcounter{equation}{0}

Abelian and non-Abelian Chern--Simons vortices have been studied in a series of relativistic or non-relativistic field models \cite{dvsc, pakh,hkp,jaw1,le1,lelw}.
These configurations  play important roles in both theoretical and experimental studies in modern physics\cite{Fro,Wil,Khf,Mat,Sok,She,Che}.
Recently the Chern--Simons  action has  received much attention in connection with  Superstring theory and M-theory \cite{sch}.
Based on the seminal work of Bagger--Lambert \cite{BL1,BL2,BL3} and Gustavsson \cite{gus} (BLG),  in  \cite{abjm}  Aharony, Bergman, Jafferis and  Maldacena
(ABJM) have  constructed  an $\mathcal{N}=6$ superconformal $U(N)\times U(N)$ Chern--Simons gauge theory
with level $(k, -k)$, coupled to four complex scalars and four fermions in the bifundamental representation.
This theory  describes the world-volume dynamics of  $N$ coincident $M2$-branes  moving in a $\mathbb{C}_4/\mathbb{Z}_k$ orbifold background in M-theory.

For both  models (BLG  and  ABJM),  it  is interesting to identify  soliton-like objects,  such as domain walls, vortices, and Q-balls,   since they correspond to various configurations of membranes.
To  this purpose, in \cite{kkkn} the authors derived  self-dual equations of BPS-type  \cite{bo,ps}  for the mass-deformed  ABJM model with gauge group $U(N) \times U(N)$.
When  $N=2$   such  BPS equations reduce  to  the Abelian Higgs  equation \cite{GST,jata,NO} whose existence theory is completely  understood \cite{jata}. While for $N>2$,
the BPS equations   give rise   to more general systems of the type describing  non-Abelian vortex  configurations,  for which   existence
and uniqueness results   were   established recently in  \cite{hy1}.
We mention also \cite{ak},  where another  set of   BPS equations were  derived  in the mass-deformed  ABJM model at both weak and strong coupling,
 in this case existence and uniqueness results  were established in \cite{liey1}.  For more recent developments  concerning the ABJM model see \cite{hkkt,kim1,bjlr,pasa,muna1,napa1,lnrm,mnrs,anru,mosu,nosht} and the references therein.

 In order to find  possible applications of the ABJM theory to specific condensed matter systems, recently a consistent  Abelian truncations of the mass-deformed ABJM model was introduced by
  Mohammed--Murugan--Nastse \cite{mmn1,mmn2}. By means of   suitable ansatz involving  fuzzy sphere matrices, the authors of \cite{mmn1,mmn2} were able
 to reduce the model into an effective  Abelian field model,  which admits a BPS reduction. However, due to the  specific  ansatz introduced in \cite{mmn1,mmn2}, the corresponding system of self-dual equations admits new features compared to those obtained
 in \cite{ak,kkkn}. In particular, we see that, to attain finite energy, the corresponding self-dual vortex configurations can be only of non-topological type. This is a surprising  new
 feature,  when  compared to  more familiar  BPS equations arising both in Abelian and non-Abelian contexts also in presence of the Chern--Simons terms, see \cite{dvsc,D1,hkp, chjde, haly, hata, HLTY, jaw1, le1, lelw, liey1, lpyjfa, pakh, tabk, yang1}
 and  references therein.  Even more importantly, for the BPS system in \cite{mmn2} we observe no energy-gap featuring along its solutions, as in fact each  (positive) energy
 level is always attained by a (constraint) energy minimizer.  The main purpose of this paper is to establish with  mathematical  rigor all such new features
 characterizing the Abelian truncation of the mass-deformed ABJM model proposed in \cite{mmn1,mmn2}. Thus, we shall prove   existence and multiplicity results
  for the  system of  BPS equations  in \cite{mmn2} in terms of the corresponding energy  and fluxes.
  As usual, to tackle such issues, we  reduce  the self-dual equations into   a planar $2\times2$ system of nonlinear elliptic equations,
  except that  now we face  an  ``indefinite '' system of elliptic PDE's,  particularly delicate to handle  by the available  analytical and variational tools.
   Indeed, so far we are aware only of the numerical study provided in \cite{mmn2}.

We shall focus to the case where all vortex points are superimposed (with assigned multiplicity) at a given point, and thus we   search for  radially symmetric solutions
about such point.  We point  out that even the radial solvability of the given system of PDE's  is non-trivial. In fact,  we lack the standard a priori estimates, so that not all   entire radial solutions satisfying non-topological
boundary conditions give rise to finite-energy configurations. However,  we manage to show that some partial integrability can be  always guaranteed. In addition, we shall  construct
a first  class of finite-energy (radial) solutions, which enjoy a physically interesting ``concentration'' property.  To this end, we use    a perturbation approach introduced in \cite{chaim} to  deal with Abelian
non-topological Chern--Simons vortices, and further developed in \cite{chta1,chta2}  for self-dual electroweak models with and without gravitational effects.
More precisely, we shall   identify   an  appropriate new   scaling  for the solution, which,  in one hand, sets our problem into a ``perturbation'' framework,
and at the same time  also helps to clarify its  specific analytic  features.  With this information in hand, and by means of a blow-up analysis,  we arrive to identify necessary and
 sufficient conditions, which characterize  the sharp range  attained by the total energy (or by  the two fluxes) along BPS-solutions.  In fact, for each value  in this range, we  construct  a 1-parameter
 family of solutions carrying such fixed amount of total  energy or flux.   This indicates that,  for the given model, the set of total fluxes covers  a planar region bounded
 by two curves emanating from a zero energy configuration. See Theorem \ref{th1} and \ref{th2} for the precise statements.  In this way, we record the  absence of ``quantization''
 effects as well as  energy gaps  along BPS solutions.

The rest of this  paper is organized as follows. In Section \ref{sec2}, (following \cite{mmn2}) we review the  ABJM model, the associated equations of motion,  the
 Abelian truncation and the derivation of the self-dual equations,  and state our  existence results.
In Section \ref{sec3}  we present necessary conditions  for  the existence of finite-energy self-dual solutions. To show their existence,
in section \ref{sec4} we introduce a new scaling, which enables  us to view  our problem as  a perturbation of a Liouville  system. So we shall apply the perturbation  approach  introduced
in \cite{chaim} to  construct a class of  radially symmetric solutions, which however do not exhaust the full range of admissible energies (or fluxes). Thus, to complete our
analysis, we provide in Section \ref{sec5} some general information about   radial solutions, and in particular we analyze   the Cauchy problem and realize that,  the finite-energy conditions
 may not be always  attained.  Therefore,  by means of  a  blow-up analysis in  Section \ref{sec6}, we shall be able to identify suitable set of  initial data  which yield to the desired  1-parameter family of finite-energy solutions.

 \section{Self-dual equations  and existence theorems}\label{sec2}
 \setcounter{equation}{0}

Let $g_{\mu\nu}={\rm diag}(+1, +1, -1)$ be the metric in the Minkowski spacetime $\mathbb{R}^{2,1}$. In the following we use (unless otherwise specified) the standard notation of summation
over repeated lower and upper indices,  and where indices are raised or lowered by using the metric in the usual way. The ABJM model \cite{abjm} is an ${\cal N}=6$ supersymmetric $U(N)\times U(N)$ Chern--Simons gauge theory with level $(k,-k)$ coupled to four  complex scalars
$C^I$ and four Dirac fermions $\psi_I$  in the fundamental representation of the $SU(4)_R$ symmetry group,  where the  R-symmetry index $I=1,...,4$.  Also,  we denote the  gauge fields for
 the two groups respectively  $A_\mu$ and $\hat A_\mu$, with Lorentz index $\mu=0,1,2$. The ABJM action  is defined as follows:
\bea
\mathcal{A}&=&\int (\mathcal{L}_{\rm CS}+\mathcal{L}_{\mbox{kin}}-V_{\mbox{ferm}}-V_0)\ud x,\label{lg1}
\eea
where  ${\cal L}_{\rm CS}$  denotes  the Chern--Simons Lagrangian  density:
\ber
\mathcal{L}_{\rm CS}&=&\frac{k}{4\pi}\epsilon^{\mu\nu\lambda}\Tr\Big(A_{\mu}\partial_{\nu}A_{\lambda}+\frac{2i}{3}A_{\mu}A_{\nu}A_{\lambda} -\hat{A}_{\mu}\partial_{\nu}\hat{A}_{\lambda}-\frac{2i}{3}\hat{A}_{\mu}\hat{A}_{\nu}\hat{A}_{\lambda}\Big);\label{lcs}
\eer
and the matter-kinetic Lagrangian density ${\cal L}_{\rm kin}$ is given by:
\ber
    \mathcal{L}_{\rm kin}&=&-\Tr\big(D_{\mu}C_{I}^{\dagger}D^{\mu}C^{I}\big)
      -i\Tr\big(\psi^{I\dagger}\gamma^{\mu}D_{\mu}\psi_{I}\big),\label{i3}
\eer
 with $\epsilon^{\mu\nu\lambda}$   the antisymmetric Levi--Civita tensor fixed by setting $\epsilon^{012}=1$,  and  $k>0$   the Chern--Simons level.
  In addition, the Dirac matrices  $\gamma^\mu$   are expressed in terms of   the standard Pauli spin matrices,
$\gamma^0=\ri\sigma^2,\gamma^1=\sigma^1,\gamma^2=\sigma^3$,
 and  the gauge-covariant derivatives are  defined by:
\ber
 D_{\mu} C^{I}&=&\partial_{\mu}C^{I}+\ri \left(A_{\mu}C^{I}-C^{I}\hat{A}_{\mu}\right),\label{i4}\\
D_\mu \psi_I &=&\partial_\mu \psi_I+\ri \left(A_\mu \psi_I- \psi_I\hat{A}_\mu\right).\label{i5}
  \eer
Furthermore, in \eqref{lg1},  $V_{\rm ferm}$   is the Yukawa-type quartic-interaction potential density given by:
\ber
  V_{\rm ferm}&=&\frac{2\pi i}{k}\Tr\Big(C_{I}^{\dagger}C^{I}\psi^{J\dagger}\psi_{J}-\psi^{J\dagger}C^{I}C_{I}^{\dagger}\psi_{J}
      -2C_{I}^{\dagger}C^{J}\psi^{I\dagger}\psi_{J}+2\psi^{J\dagger}C^{I}C_{J}^{\dagger}\psi_{J}\nm\\
&&+\epsilon^{IJKL} C_{I}^{\dagger}\psi_{J}C_{K}^{\dagger}\psi_{L}-\epsilon_{IJKL}C^{I}\psi^{J\dagger}C^{K}\psi^{L\dagger}\Big);\label{i6}
   \eer
 and  $V_0$ is the sextic scalar potential:
 \ber
 V_0&=& \frac{4\pi^2}{3k^2}\Tr \Big(C^{I}C_{I}^{\dagger}C^{J}C_{J}^{\dagger}C^{K}C_{K}^{\dagger}
+C_{I}^{\dagger}C^{I}C_{J}^{\dagger}C^{J}C_{K}^{\dagger}C^{K}+4C^{I}C_{J}^{\dagger}C^{K}C_{I}^{\dagger}C^{J}C_{K}^{\dagger}\nm\\
&& -6C^{I}C_{J}^{\dagger}C^{J}C_{I}^{\dagger}C^{K}C_{K}^{\dagger}C^{K}\Big),\label{i7}
   \eer
   where again:  $\epsilon^{1234}=1$,   and the antisymmetric tensor $\epsilon^{IJKL}$ takes values $\pm1$ according to the whether the set of indices  $\{IJKL\}$ is an even or odd permutation of $\{1234\}$
   and it is set to  zero otherwise.

   The ABJM action exhibits  an $SU(4)\times U(1)$ R-symmetry associated with the   $\mathcal {N}=6$ supersymmetric  transformation \cite{abjm,kkkn}.
 This theory has a maximally supersymmetric massive deformation with mass parameter $\sigma$ \cite{grvv,tera},
 which breaks the R-symmetry down to $SU(2) \times SU(2)\times U(1)_{A}\times U(1)_{\hat{A}}\times \mathbb{Z}_{2}$  by  decomposing the scalars as follows:
\be
C^{I}=(Q^{\alpha},R^{\alpha}), \quad \alpha=1, 2. \label{i9}
\ee

In this way the $\mathbb{Z}_2$-action swaps the matter fields $Q^\alpha$ and $R^\alpha$, the $SU(2)$-action is taken individually on $Q^\alpha$ and $R^\alpha$
respectively, while the $U(1)$ symmetry rotates $Q^\alpha$ on phase $+1$ and $R^\alpha$ by a phase $-1$.
The bosonic part of the mass-deformed Lagrangian takes the form \cite{mmn1,mmn2}:
\ber
\mathcal{L}_{\rm Bosonic}&=&\frac{k}{4\pi}\epsilon^{\mu \nu \lambda}\Tr\left(A_{\mu}\partial_{\nu}
  A_{\lambda}+\frac{2\ri}{3}A_{\mu}A_{\nu}A_{\lambda}-\hat{A}_{\mu}\partial_{\nu}\hat{A}_{\lambda}
  -\frac{2\ri}{3}\hat{A}_{\mu}\hat{A}_{\nu}\hat{A}_{\lambda}\right)\nonumber\\
  &-&\Tr\big(|D^{\mu}Q^{\alpha}|^2\big)-\Tr\big(|D^{\mu}R^{\alpha}|^2\big)-V, \label{i10}
\eer
in which  the potential $V$ takes the following sextic form:
\ber
     V&=&\Tr\left(|M^{\alpha}|^2+|N^{\alpha}|^2\right),\label{i11}
\eer
with
 \ber
 M^{\alpha}=\sigma Q^{\alpha}+\frac{2\pi}{k}\Big(2Q^{[\alpha}Q^{\dagger}_{\beta}Q^{\beta]}+R^{\beta}R^{\dagger}_{\beta}Q^{\alpha}-Q^{\alpha}R^{\dagger}_{\beta}R^{\beta}
   +2Q^{\beta}R^{\dagger}_{\beta}R^{\alpha}-2R^{\alpha}R^{\dagger}_{\beta}Q^{\beta}\Big), \label{i12}\\
   N^{\alpha} =-\sigma R^{\alpha}+\frac{2\pi}{k}\Big(2R^{[\alpha}R^{\dagger}_{\beta}R^{\beta]}+Q^{\beta}Q^{\dagger}_{\beta}R^{\alpha}-R^{\alpha}Q^{\dagger}_{\beta}Q^{\beta}
   +2R^{\beta}Q^{\dagger}_{\beta}Q^{\alpha}-2Q^{\alpha}Q^{\dagger}_{\beta}R^{\beta}\Big), \label{i13}
 \eer
 and where  the notation  $Y^{[\alpha}Y^\dagger_\gamma Y^{\beta]}\equiv\left(Y^\alpha Y^\dagger_\gamma Y^\beta-Y^\beta Y^\dagger_\gamma Y^\alpha\right)$ is observed.

The  Euler--Lagrange equations of  the bosonic Lagrangian \eqref{i10}  are given by \cite{mmn1,mmn2}:
\ber
 \frac{k}{2\pi}\epsilon^{\mu\nu\lambda}F_{\nu\lambda}&=&J^{\mu},\label{i14}\\
  \frac{k}{2\pi}\epsilon^{\mu\nu\lambda}\hat{F}_{\nu\lambda}&=&\hat{J}^{\mu},\label{i15}\\
  D_{\mu}D^{\mu}Q^{\alpha}&=&\frac{\partial V}{\partial Q^{\dagger}_{\alpha}}, \label{i16}\\
  D_{\mu}D^{\mu}R^{\alpha}&=&\frac{\partial V}{\partial R^{\dagger}_{\alpha}},\label{i17}
\eer
with field strengths:
\ber
F_{\mu \nu}=\partial_{\mu}A_{\nu}-\partial_{\nu}A_{\mu}+\ri[A_{\mu}, A_{\nu}], \quad \hat{F}_{\mu \nu}=\partial_{\mu}\hat{A}_{\nu}-\partial_{\nu}\hat{A}_{\mu}+\ri[\hat{A}_{\mu}, \hat{A}_{\nu}], \label{i18}
 \eer
and   the two currents $J^{\mu}$ and $\hat{J}^{\mu}$  defined by
\ber
  J^{\mu}&\equiv&\ri\left(Q^{\alpha}(D^{\mu}Q^{\alpha})^{\dagger}-(D^{\mu}Q^{\alpha})Q_{\alpha}^{\dagger}
  +R^{\alpha}(D^{\mu}R^{\alpha})^{\dagger}-(D^{\mu}R^{\alpha})R_{\alpha}^{\dagger}\right), \label{i19}\\
  \hat{J}^{\mu}&\equiv&-\ri\left(Q_{\alpha}^{\dagger}(D^{\mu}Q^{\alpha})-(D^{\mu}Q^{\alpha})^{\dagger}Q^{\alpha}
  +R_{\alpha}^{\dagger}(D^{\mu}R^{\alpha})-(D^{\mu}R^{\alpha})^{\dagger}R^{\alpha} \right),\label{i20}
\eer
which are  covariantly conserved: $\nabla_{\mu}J^{\mu}=\nabla_{\mu}\hat J^\mu=0$. In addition,  there are two  ordinarily conserved Abelian currents $j^{\mu}$  and $\hat{j}^{\mu}$:
\ber
   j^{\mu}&\equiv&\ri\Tr\left(Q^{\alpha}(D^{\mu}Q^{\alpha})^{\dagger}-(D^{\mu}Q^{\alpha})Q_{\alpha}^{\dagger}+R^{\alpha}(D^{\mu}R^{\alpha})^{\dagger}-(D^{\mu}R^{\alpha})R_{\alpha}^{\dagger}\right), \label{i21}\\
  \hat{j}^{\mu}&\equiv&-\ri\Tr\left(Q_{\alpha}^{\dagger}(D^{\mu}Q^{\alpha})-(D^{\mu}Q^{\alpha})^{\dagger}Q^{\alpha}+R_{\alpha}^{\dagger}(D^{\mu}R^{\alpha})-(D^{\mu}R^{\alpha})^{\dagger}R^{\alpha}\right),\label{i22}
\eer
with respect to the global $U(1)_{A}$ and $U(1)_{\hat{A}}$ invariances.

In view of the complicated structures of the equations \eqref{i14}--\eqref{i17}, one is prompted  to consider some simplification.
More precisely,  by introducing the  real-valued vector fields $a_{\mu}^{(1)}$ and $a_{\mu}^{(2)}$,  and the complex-valued scalar fields $\phi_{\alpha}$ and  $\chi_{\alpha}$,
 Mohammed--Murugan--Nastase \cite{mmn1,mmn2}  proposed  the following Abelianization    ansatz for the Chern--Simons fields and the matter scalars:
\ber
A_{\mu}&=&a^{(2)}_{\mu}G^{1}G_{1}^{\dagger}+a^{(1)}_{\mu}G^{2}G_{2}^{\dagger},\label{i24}\\
\hat{A}_{\mu}&=&a^{(2)}_{\mu}G_{1}^{\dagger}G^{1}+a^{(1)}_{\mu}G_{2}^{\dagger}G^{2},\label{i25}\\
Q^{\alpha}&=&\phi_{\alpha}G^{\alpha},\label{i26}\\
R^{\alpha}&=&\chi_{\alpha}G^{\alpha},\label{i27}
\eer
$\alpha=1,2$ and $\mu=0,1,2$ (no summation over $\alpha$ in \eqref{i26}--\eqref{i27}),  with  $G^\alpha$  defined  as follows:
\ber
(G^1)_{m,n}=\sqrt {m-1}\delta_{m, n},\quad
(G^2)_{m,n} =\sqrt {N-m}\delta_{m+1, n}, \label{i28}\\
 (G_1^{\dagger})_{m,n}=\sqrt{m-1} ~\delta_{m, n}, \quad
(G_2^{\dagger})_{m,n}=\sqrt{N-n} ~\delta_{n+1, m},\label{i29}
\eer
$m,n\in\{1,\dots, N\}$.  In particular, the mass-deformed theory admits the following  ground state of the fuzzy type \cite{grvv,tera}:
\be
R^\alpha=cG^\alpha,\qquad  Q^\alpha=0\qquad \text{and}\qquad  Q^\dagger_\alpha=cG^\alpha,\qquad R^\alpha=0,\label{i30}
\ee
with
\ber
c\equiv\sqrt{\frac{\sigma k}{2\pi}}.\label{i31}
\eer

Actually, the ground state \eqref{i30} corresponds to a fuzzy 2-sphere \cite{npr1, np1}.
Note also that the matrices $G^\alpha$, bifundamental under $U(N)\times U(N)$, satisfy:
\be
G^\alpha=G^\alpha G^\dagger_\beta G^\beta-G^\beta G^\dagger_\beta G^\alpha,\quad \alpha,\beta=1,2, \label{i32}
\ee
(again no summation on repeated indices  in \eqref{i32}).

The above ansatz leads to  a consistent Abelian  truncation of the ABJM Lagrangian \cite{mmn1,mmn2},
\ber
\mathcal{L}_{\rm td}=-\frac{N(N-1)}{2}\Bigg[\frac{k}{4\pi}\epsilon^{\mu \nu \lambda}\left(a^{(2)}_{\mu}f^{(1)}_{\nu \lambda}+a^{(1)}_{\mu}f^{(2)}_{\nu \lambda}\right)
+|D_{\mu}\phi_{i}|^{2}+|D_{\mu}\chi_{i}|^{2}+U(|\phi_{i}|,|\chi_{i}|)\Bigg], \label{i33}
\eer
 with the following reduced potential  $U\equiv 2V/N(N-1)$:
\bea
U&=&\frac{4\pi^{2}}{k^2}\Big[\left(|\phi_{1}|^{2}+|\chi_{1}|^{2}\right)\left(|\chi_{2}|^{2}-|\phi_{2}|^{2}-c^{2}\right)^{2} +4|\phi_{1}|^{2}|\phi_{2}|^{2}\left(|\chi_{1}|^{2}+|\chi_{2}|^{2}\right)\nm\\
&&+(|\phi_{2}|^{2}+|\chi_{2}|^{2})\left(|\chi_{1}|^{2}-|\phi_{1}|^{2}-c^{2}\right)^{2}+4|\chi_{1}|^{2}|\chi_{2}|^{2}\left(|\phi_{1}|^{2}+|\phi_{2}|^{2}\right)\Big]. \label{i34}
\label{abelianpot}
\eea
Clearly,  the induced   Abelian gauge covariant derivatives are given by:
\ber
D_{\mu}\phi_{i}=\partial_{\mu}\phi_{i}-\ri a_{\mu}^{(i)}\phi_{i},\quad D_{\mu}\chi_{i}=\partial_{\mu}\chi_{i}-\ri a_{\mu}^{(i)}\chi_{i}, \label{i35}
\eer
with  the  induced Abelian field strengths:
\ber
 f_{\mu\nu}^{(i)}=\partial_\mu a^{(i)}_\nu-\partial_\nu a^{(i)}_\mu, \quad i=1, 2. \label{i36}
\eer

The above Abelianization involves four complex-valued scalars, and  this  general situation  can be simplified further by setting: $\chi_1=\phi_2=0$.
Then, with abuse of notation, where one denotes $\phi_2$ in place of $\chi_2$, in  \cite{mmn1,mmn2}  Mohammed, Murugan and Nastase  obtained
the following effective  Lagrangian:
\ber
\mathcal{L}_{\rm td}&=&-\frac{N(N-1)}{2}\left[\frac{k}{4\pi}\epsilon^{\mu \nu \lambda}\left(a^{(2)}_{\mu}f^{(1)}_{\nu \lambda}+a^{(1)}_{\mu}f^{(2)}_{\nu \lambda}\right)
+|D_{\mu}\phi_{i}|^{2} +U(|\phi_1|, |\phi_2|)\right],\label{i37}
\eer
with
\ber
U(|\phi_1|, |\phi_2|)&=&\frac{4\pi^2}{k^2}\left[|\phi_1|^2(|\phi_2|^2-c^2)^2+|\phi_2|^2(|\phi_1|^2+c^2)^2\right],\label{i38}
\eer
and $c$  given in \eqref{i31}.

 It  can be easily checked that the    equations of motion  for the truncated  Lagrangian \eqref{i37} take the form:
\ber
 \frac{k}{2\pi}\epsilon^{\mu\nu\alpha}f^{(1)}_{\nu\alpha}&=&\ri\left(\phi_2\overline{D^\mu\phi_2}-\overline{\phi_2}D^\mu\phi_2\right),\label{i39}\\
  \frac{k}{2\pi}\epsilon^{\mu\nu\alpha}f^{(2)}_{\nu\alpha}&=&\ri\left(\phi_1\overline{D^\mu\phi_1}-\overline{\phi_1}D^\mu\phi_1\right),\label{i40}\\
  D_\mu D^\mu\phi_1&=&\frac{4\pi^2}{k^2}\left[\left(|\phi_2|^2-c^2\right)^2+2|\phi_2|^2\left(|\phi_1|^2+c^2\right)\right]\phi_1,\label{i41}\\
   D_\mu D^\mu\phi_2&=&\frac{4\pi^2}{k^2}\left[2|\phi_1|^2\left(|\phi_2|^2-c^2\right)+\left(|\phi_1|^2+c^2\right)^2\right]\phi_2.\label{i42}
\eer

We are interested  to obtain  static solutions of  \eqref{i39}--\eqref{i42}.
To this purpose, we use the  Gauss laws of the Lagrangian \eqref{i37} (expressed by the $\mu=0$  component of  \eqref{i39}--\eqref{i40}),  given as
follows:
\ber
\frac{k}{2\pi}f_{12}^{(1)}=a_0^{(2)}|\phi_2|^2, \label{gs1}\\
\frac{k}{2\pi}f_{12}^{(2)}=a_0^{(1)}|\phi_1|^2, \label{gs2}
\eer     together with the following well-known identities:
\ber
 |D_1\phi_i|^2+ |D_2\phi_i|^2=|D_1\phi_i\pm\ri D_2\phi_i|^2\pm\ri\left(\partial_1\big[\phi_i\overline{D_2\phi_i}\big]-\partial_2\big[\phi_i\overline{D_1\phi_i}\big]\right)\pm f^{(i)}_{12}|\phi_i|^2,\,\, i=1, 2.\, \label{i43}
\eer
As a consequence, we can  write the static energy density relative to \eqref{i37} as follows:
 \ber
 \mathcal{E}&=&-\mathcal{L}_{\rm td}=\frac{N{(N-1)}}{2}\Bigg\{\left(a_0^{(1)}\right)^2|\phi_1|^2+\left(a_0^{(2)}\right)^2|\phi_2|^2+\sum\limits_{i=1}^2\left(|D_i\phi_1|^2+|D_i\phi_2|^2\right)\nm\\
 &&+\frac{4\pi^2}{k^2}\Big[|\phi_1|^2(|\phi_2|^2-c^2)^2+|\phi_2|^2(|\phi_1|^2+c^2)^2\Big]\Bigg\}\nm\\
 &=&\frac{N{(N-1)}}{2}\Bigg\{\frac{\left(\frac{k}{2\pi}f_{12}^{(2)}\right)^2}{4|\phi_1|^2}+\frac{\left(\frac{k}{2\pi}f_{12}^{(1)}\right)^2}{4|\phi_2|^2}
 +\sum\limits_{i=1}^2\left(|D_i\phi_1|^2+|D_i\phi_2|^2\right)\nm\\
 &&+\frac{4\pi^2}{k^2}\Big[|\phi_1|^2(|\phi_2|^2-c^2)^2+|\phi_2|^2(|\phi_1|^2+c^2)^2\Big]\Bigg\}\nm\\
 &=&\frac{N{(N-1)}}{2}\Bigg\{\left(\frac{\frac{k}{2\pi}f_{12}^{(2)}}{2|\phi_1|}\pm\frac{2\pi}{k}|\phi_1|(|\phi_2|^2-c^2)\right)^2+\left(\frac{\frac{k}{2\pi}f_{12}^{(1)}}{2|\phi_2|}\pm\frac{2\pi}{k}|\phi_2|(|\phi_1|^2+c^2)\right)^2\nm\\
 &&+|D_1\phi_1\pm\ri D_2\phi_1|^2+|D_1\phi_2\pm\ri D_2\phi_2|^2\pm\ri\left(\partial_1\big[\phi_1\overline{D_2\phi_1}\big]-\partial_2\big[\phi_1\overline{D_1\phi_1}\big]\right)\nm\\
 &&\pm\ri\left(\partial_1\big[\phi_2\overline{D_2\phi_2}\big]-\partial_2\big[\phi_2\overline{D_1\phi_2}\big]\right)\pm c^2\left(f_{12}^{(2)}-f_{12}^{(1)}\right)\Bigg\}. \label{i44}
 \eer
In this way  we see that,  when we can neglect the pure divergence terms, then  the corresponding total  energy  admits the following  lower bound:
 \ber
  E=\int \mathcal{E} \ud x \ge \pm\frac{N{(N-1)}}{2}c^2\int \left(f_{12}^{(2)}-f_{12}^{(1)}\right)\ud x, \label{i45}
 \eer
 which is  attained by solutions of   the following  first-order  self-dual equations:
 \ber
  &&D_1\phi_1\pm\ri D_2\phi_1=0, \label{i46}\\
   &&D_1\phi_2\pm\ri D_2\phi_2=0, \label{i47}\\
   &&f_{12}^{(1)}\pm \frac{8\pi^2}{k^2} |\phi_2|^2\left(|\phi_1|^2+c^2\right)=0,\label{i48}\\
   &&f_{12}^{(2)}\pm\frac{8\pi^2}{k^2}|\phi_1|^2\left(|\phi_2|^2-c^2\right)=0, \label{i49}
 \eer
 which need to be supplemented by the Gauss-law constraints \eqref{gs1}--\eqref{gs2}.

 As a consequence, to obtain a self-dual field configurations, it suffices to solve for $(a_1^{(i)}, a_2^{(i)}, \phi_i) $ $(i=1,2)$ the self-dual equations \eqref{i46}--\eqref{i49}
 and then use \eqref{gs1}--\eqref{gs2} in order to obtain the remaining component $a_0^{(i)}$ $ (i=1,2)$. Such   configurations identify a special class of
 solutions for the more general field equations \eqref{i39}--\eqref{i42}.

 The BPS-type equations \eqref{i46}--\eqref{i49} were first  derived  in \cite{mmn2}, where the authors also  provided a discussion and some
 numerical evidence about their solvability.  Our goal in  this paper is to establish with mathematical rigor existence and multiplicity results about solutions
 of   \eqref{i46}--\eqref{i49},  also in terms of the corresponding  energy levels.

 Actually, compared to  other well known  self-dual equations of BPS type governing Abelian and non-Abelian vortex configurations,  we shall  see that
 the set of equations  \eqref{i46}--\eqref{i49} can  support only planar  non-topological type solutions.  In fact,  the more classical  topological solutions (typical of  the
 Abelian Higgs model) are not admissible here. So, in this case,  no special ``quantization'' of energy can be identified in terms of energy minimizer over appropriate
 topological classes. This fact also excludes the presence of an energy gap along solutions of \eqref{i46}--\eqref{i49}.
   On the contrary, we show that any (positive) energy level can be attained together with the equal sign in \eqref{i45}.

In addition, we point out that,  for planar solutions of  \eqref{i46}--\eqref{i49}, the usual non-topological boundary condition:
 \ber
 |\phi_i|^2\to 0\quad \text{as}\quad |x|\to \infty,\quad i=1, 2,\label{i50}
 \eer
 is  necessary  but not sufficient   to accomplish the finite-energy condition.

 In this paper we shall be able to  handle such a new situation  in case the complex scalar fields admit the same unique zero with possible different multiplicities.
Recall that,  on the basis of equations  \eqref{i46} and \eqref{i47},   we know  that the zeros of the complex scalars $\phi_1$ and $\phi_2$  are discrete and
 with integral  multiplicities (see \cite{JT}). Furthermore, in a self-dual regime, such zeros do not interact, and so they can be arbitrarily prescribed.
 Thus, we focus to  the case where both $\phi_1$ and $\phi_2$ admit the same  unique zero at a point $p$ with multiplicities $N_1$ and $N_2$, respectively.
This fact will allow us to look for solutions radially symmetric about such point.

Since  solutions of   equations \eqref{i46}--\eqref{i49} with the lower  and upper sign are related via  the  relation   $(a^{(i)}_1, a^{(i)}_2, \phi_i)\leftrightarrow (-a^{(i)}_1, -a^{(i)}_2, \overline{\phi_i})$ $(i=1, 2)$,
in what follows  we only  deal with  the case  where we specify the  lower sign in \eqref{i46}--\eqref{i49}.

 We  show that in such case   the total fluxes:
 \be
 \Phi_i=\int_{\mathbb{R}^2}f^{(i)}_{12}\ud x,\quad i=1,2\label{fl}
 \ee
 must satisfy:
 \be
 \frac{1}{2\pi} \Phi_2\in(N_2+1,\,\, N_1+N_2+1) \quad \text{and}\quad \frac{1}{2\pi} \Phi_1\in(N_1+N_2+1,\,\, \infty). \label{fl1}
 \ee

 However, such set of conditions are not sufficient, in the sense that we cannot assign  independently the values of $\Phi_1$ and $\Phi_2$, as prescribed by \eqref{fl1}.

  We shall discuss below to what extent the necessary conditions \eqref{fl1} turn into sufficient ones. But first to check that $\Phi_1$ and $\Phi_2$ influence each other,
   we state the following result which describes a very special class of field configurations,  solutions of \eqref{i46}--\eqref{i49},   characterized by a physically relevant  ``concentration'' property.

 \begin{theorem}\label{th0}
  Let $k>0$, $\sigma>0$, $p\in \mathbb{R}^2$ and $N_1,N_2\in \mathbb{N}$. There exists $\vep_0>0$ such that $\forall\,\vep\in(0, \vep_0)$ the system \eqref{i46}--\eqref{i49}
  admits a field configuration:  $(a^{(i),\vep}_1, a^{(i),\vep}_2, \phi_i^\vep)$ $(i=1,2)$  where  $p$ is the only zero of $\phi_i^\vep$ with multiplicity $N_i$ $(i=1,2)$, and
   \ber
   |\phi_i^\vep|^2={\rm O}\left(|x|^{-2\beta_i^\vep}\right), \, |D_j\phi_i^\vep|^2={\rm O}\left(|x|^{-2(\beta_i^\vep+1)}\right), \,f^{(i),\vep}_{12}={\rm O}\left(|x|^{-2\beta_i^\vep}\right),\,\text{as}\, |x|\to\infty \label{2.53a}
     \eer
  with $\beta_1^\vep=N_2+1+{\rm o}(1)$ and $\beta_2^\vep=N_1+1+{\rm o}(1)$ as $\vep\to0^+$.

  Furthermore, as   $\vep\to0^+$:
   \ber
    &&\frac{1}{2\pi}f^{(i), \vep}_{12}\rightharpoonup (N_1+N_2+1)\delta_p\quad \text{weakly in the sense of measure in }\,\,\mathbb{R}^2\,\,(i=1,2),\label{2.54}\\
    &&f^{(1),\vep}_{12}-f^{(2),\vep}_{12}\to 0 \quad \text{in }\quad L^1(\mathbb{R}^2). \label{2.55}
   \eer
   In particular, as   $\vep\to0^+$:
   \ber
    &&\frac{1}{2\pi} \Phi_{i,\vep}=\frac{1}{2\pi}\itr f^{(i), \vep}_{12}\ud x\to N_1+N_2+1 \,\,(i=1,2), \label{2.56}\\
    &&E_\vep= \frac{N(N-1)\sigma k}{2\pi}\itr \left(f^{(1),\vep}_{12}-f^{(2),\vep}_{12}\right)\ud x\to 0. \label{2.57}
   \eer

 \end{theorem}
\hfill $\square$

Notice that  the configurations  constructed in Theorem \ref{th0}  cover only a very tiny part of full range of values specified in \eqref{fl1}, however,
they confirm the fact that the values of the fluxes $\Phi_1$ and $\Phi_2$ are tied to each other.   This suggests that we are free to prescribe only one of the
two fluxes,  according to \eqref{fl1}.  We prove this fact in the following:

 \begin{theorem}\label{th1}
Let $k>0$ and $\sigma>0$. For any point $p\in \mathbb{R}^2$ and $N_1,N_2\in \mathbb{N}$,  the condition:
 \be
 \frac{1}{2\pi} \Phi_2\in(N_2+1,\,\, N_1+N_2+1) \quad\bigg( \text{or}\quad \frac{1}{2\pi} \Phi_1\in(N_1+N_2+1,\,\, \infty)\bigg) \label{fl2}
 \ee
is necessary and sufficient for the existence of a finite-energy field configuration $(a^{(i)}_1, a^{(i)}_2, \phi_i)$ $(i=1,2)$ solutions of \eqref{i46}--\eqref{i49} in $\mathbb{R}^2$,
with $\phi_i$ vanishing exactly at $p$ with multiplicities $N_i (i=1,2)$.

Furthermore, the following  asymptotic estimates hold:
    \ber
   |\phi_i|^2={\rm O}\left(|x|^{-2\beta_i}\right), \, |D_j\phi_i|^2={\rm O}\left(|x|^{-2(\beta_i+1)}\right), \,f^{(i)}_{12}={\rm O}\left(|x|^{-2\beta_i}\right),\,\text{as}\, |x|\to\infty \label{2.54'}
     \eer
 with $\beta_i\equiv\frac{1}{2\pi} \Phi_i-N_i$  ($i=1, 2$) satisfying:  $\beta_1\in(N_2+1, \,\,\infty)$ and  $\beta_2\in(1,\,\, N_1+1)$, and more specifically,
  \ber
   &&0<(\beta_1-1)(\beta_2-1)<(N_1+1)(N_2+1),\label{s68a}\\
    &&\beta_2(\beta_1-1)>N_2(N_1+1), \label{s69'}\\
    &&0<\beta_1(\beta_2-1)<N_1(N_2+1).  \label{s70'}
  \eer
 \end{theorem}
%\hfill $\square$

 Theorem \ref{th0} and  Theorem \ref{th1} also   indicate that the set of admissible fluxes $(\Phi_1, \Phi_2)$ describes a region in  the plane  bounded by two curves parameterized by $\Phi_1$ (or $\Phi_2$),
 emanating form a zero-energy configuration.

  The following multiplicity  result substantiates this fact.

 \begin{theorem} \label{th2}
 Let $k>0$, $\sigma>0$, $p\in \mathbb{R}^2$ and $N_1,N_2\in \mathbb{N}$. For any $\gamma_2\in(N_2+1,\, N_1+N_2+1)$ (or $\gamma_1>N_1+1$), there exits a 1-parameter family of finite-energy solutions
 of  \eqref{i46}--\eqref{i49} in $\mathbb{R}^2$ such that the complex scalar fields $\phi_i$  vanish only  at $p$ with multiplicity $N_1$ and $N_2$ respectively,  and
 \ber
 \gamma_2= \frac{1}{2\pi} \Phi_2 \quad (\text{or}\,\,  \gamma_1= \frac{1}{2\pi} \Phi_1). \label{fl2'}
 \eer

  Equivalently, for any assigned energy level $E>0$, there exists a $1$-parameter family of solutions as described above, which all carry total energy $E$, namely:
  \ber
  E= \frac{N(N-1)\sigma k}{4\pi}\int_{\mathbb{R}^2}\left(f_{12}^{(1)}-f_{12}^{(2)}\right)\ud x.
  \eer
  For all such solutions, the asymptotic estimates \eqref{2.54'} hold together with \eqref{s68a}--\eqref{s70'}.
 \end{theorem}
\hfill $\square$

 To establish the above results, we reduce   \eqref{i46}--\eqref{i49} into a $2\times2$ system of second-order elliptic equations.
 To this purpose,   let the zeros of $\phi_i$ have multiplicity $N_i\in \mathbb{N}$ $(i=1, 2)$ and  be  given  as follows:
 \ber
  p_{i1,}, \cdots,   p_{iN_i}, \quad i=1, 2, \label{2.66}
 \eer
  each repeated according to its multiplicity.

Set
  \ber
   u=\ln|\phi_1|^2-\ln\frac{\sigma k}{2\pi}, \quad v=\ln|\phi_2|^2-\ln\frac{\sigma k}{2\pi}.
  \eer
By means of  the  relations (see e.g. \cite{yang1}):
    \ber
\ds&& \phi_1=\sqrt{\frac{\sigma k}{2\pi}}\exp\left(\frac12u(x)+\ri\sum\limits_{s=1}^{N_1}\arg(x-p_{1s})\right), \label{re1}\\
\ds&& \phi_2=\sqrt{\frac{\sigma k}{2\pi}}\exp\left(\frac12v(x)+\ri\sum\limits_{s=1}^{N_2}\arg(x-p_{2s})\right), \label{re1'} \\
\ds&& a_1^{(1)}(x)=-{\rm Re}\{\ri \bar{\pa}\ln\phi_1(x)\},\quad a_2^{(1)}(x)=-{\rm Im}\{\ri \bar{\pa}\ln\phi_1(x)\},\label{re2}\\
\ds&& a_1^{(2)}(x)=-{\rm Re}\{\ri \bar{\pa}\ln\phi_2(x)\},\quad a_2^{(2)}(x)=-{\rm Im}\{\ri \bar{\pa}\ln\phi_2(x)\},\label{re2'}\\
\ds&&|\phi_1|^2=\frac{\sigma k}{2\pi}\re^u, \quad |\phi_2|^2=\frac{\sigma k}{2\pi}\re^v,\label{re3'}\\
\ds&&|D_j\phi_1|^2=\frac{\sigma  k}{4\pi}\re^u|\nabla u|^2, \quad |D_j\phi_2|^2=\frac{\sigma  k}{4\pi}\re^v|\nabla v|^2,\label{re3}
\eer
with $ \bar{\partial} =\frac12(\partial_1+\ri \partial_2)$, it is standard to see that,
\ber
 -\Delta u=2f_{12}^{(1)}-4\pi\sum_{s=1}^{N_1}\delta_{p_{1s}} \quad \text{and}\quad  -\Delta v=2f_{12}^{(2)}-4\pi\sum_{s=1}^{N_2}\delta_{p_{2s}}, \label{ill}
\eer
where  $\delta_p$  denotes the   Dirac measure concentrated at $p\in \mathbb{R}^2$.

Therefore,  we  find that,  problem \eqref{i46}--\eqref{i49}  can be equivalently   reformulated    in terms of the unknown $(u, v)$ satisfying  the following system of
 nonlinear elliptic  equations:
  \begin{equation}\label{i64}
\begin{cases}
 \ds-\Delta u=\lambda\re^{v}\left(\re^{u}+1\right)-4\pi\sum_{s=1}^{N_1}\delta_{p_{1s}}\quad \text{in}\quad \mathbb{R}^2,  \\
     \ds-\Delta v=\lambda\re^{u}\left(\re^{v}-1\right)-4\pi\sum_{s=1}^{N_2}\delta_{p_{2s}}\quad \text{in}\quad \mathbb{R}^2,
\end{cases}
\end{equation}
with
 \ber
  \lambda\equiv 4\sigma^2>0.  \label{i66}
 \eer

 In other words, any solution of \eqref{i46}--\eqref{i49} with the zeros of $\phi_i$ given by \eqref{2.66} can be expressed in terms of
 \eqref{re1}--\eqref{re2'} with $(u, v)$ satisfying \eqref{i64}. Viceversa, from  every solution of \eqref{i64} we can recover a
 field configuration, $(a_1^{(i)}, a_2^{(i)}, \phi_i)$ (i=1, 2),  solution  of  \eqref{i46}--\eqref{i49},  via the relations  \eqref{re1}--\eqref{re2'},
 such that  the zeros of $\phi_i$ $(i=1, 2)$ are  given exactly by \eqref{2.66}.

In view of \eqref{i48}, \eqref{i49},  \eqref{fl} and \eqref{re1}-\eqref{re2'}, to  obtain a planar self-dual field configuration which carries finite total energy and fluxes, we need to solve \eqref{i64} over $\mathbb{R}^2$ under the integrability conditions:
   \ber
   \re^u,\, \re^v,\, \re^{u+v}\in L^1(\mathbb{R}^2).  \label{itg}
 \eer

 In particular, \eqref{itg} implies that we must satisfy the non-topological boundary conditions:
 \ber
 u\to-\infty,\quad v\to-\infty\quad \text{as}\quad |x|\to \infty.\label{ntp}
 \eer
 However, we anticipate that, contrary to other self-dual equations, for solutions of \eqref{i64} the  boundary conditions \eqref{ntp} do  not
 always  imply \eqref{itg}.  Furthermore, under the standard transformation (natural for Liouville systems):
 \be
  u(x)\to u\big(\frac x\vep\big)+2\ln\frac1\vep, \quad  v(x)\to v\big(\frac x\vep\big)+2\ln\frac1\vep,\label{2.78}
 \ee
 with $\vep>0$ small,  problem \eqref{i64}  transforms into a ``perturbation'' of an indefinite Liouville system, which admits no solutions
 under the integrability condition \eqref{itg}. In other words, under \eqref{itg},  for the given system  we have no hope to find solutions ``bifurcating'' from
 Liouville solutions.

 On the contrary, by changing the point of view,  we find that,  it is the ``interacting'' term
 $\re^{u+v}$ to play the role of ``leading'' term within a (more appropriate) perturbation approach.
  So, under a suitable scaling, it will be  the term $\re^{u+v}$  to    behave  as a Liouville solution.   Interestingly, such an understanding will help us also to find solutions away from the perturbation regime.

To avoid additional technicalities in carrying out such program, we will focus to  the case where all the zeros of the complex scalars are superimposed at a point, which
(after a translation), we may take to be the origin.  Namely, we set
 \[p_{i1}=\cdots=p_{iN_i}=0,\quad i=1, 2.\]

 Moreover, by rescaling:
  \ber
   u(x)\to u\big(\frac{x}{\sqrt\lambda}\big)\quad\text{and} \quad v(x)\to v\big(\frac{x}{\sqrt\lambda}\big), \label{2.78a}
   \eer
  we can further assume that $\lambda=1$.  So, we need to study the solvability of   the system:
  \ber\left\{\ba{lll}
  \ds -\Delta u&=&\re^{v}(\re^{u}+1)-4\pi N_1\delta_0\quad \text{in}\quad \mathbb{R}^2,\\[2mm]
   \ds-\Delta v&=&\re^{u}(\re^{v}-1)-4\pi N_2\delta_0\quad \text{in}\quad \mathbb{R}^2,
   \ea\right.\label{ss1}
 \eer
 with $N_1,N_2\in \mathbb{N}$, under the integrability  conditions specified in   \eqref{itg}.

  By analyzing entire radial solutions of \eqref{ss1} about the origin,  we check that, even when they realized the boundary condition \eqref{ntp},  they may fail to satisfy \eqref{itg}.  Therefore, our main task will be to identify suitable set of
  initial conditions which  yield to (radial) solutions of \eqref{ss1} satisfying \eqref{itg}, see Theorem \ref{thfr} for details.  More precisely, we shall analyze their multiplicity
  along any  assigned energy level  or at fixed fluxes.

  \section {Necessary conditions for solvability}\setcounter{equation}{0} \label{sec3}

  In this section we derive necessary conditions for the solvability of the problem:
    \ber\left\{\ba{lll} \label{ss11}
  \ds -\Delta u=\re^{v}\left(\re^{u}+1\right)-4\pi N_1\delta_0\quad \text{in}\quad \mathbb{R}^2,\\[2mm]
   \ds-\Delta v=\re^{u}\left(\re^{v}-1\right)-4\pi N_2\delta_0\quad \text{in}\quad \mathbb{R}^2,\\[2mm]
   \re^u,\,\re^v,\,\re^{u+v} \in L^1(\mathbb{R}^2).
   \ea\right.
 \eer

  To this purpose we recall that, by an approach introduced by Chen--Li \cite{CL1,CL2} for Liouville equations (based on potential estimates),  it is possible to  show that every solution $(u, v)$ of \eqref{ss11} must admit the following asymptotic behavior at $\infty$:
    \ber\left\{\ba{lll}
 \ds u(x)=-\left(\frac{1}{2\pi}\int_{\mathbb{R}^2}\re^v\left(\re^u+1\right)\ud x-2N_1\right)\ln|x|+{\rm O}(1)\quad \text{as}\quad |x|\to \infty,\\
 \ds v(x)=-\left(\frac{1}{2\pi}\int_{\mathbb{R}^2}\re^u\left(\re^v-1\right)\ud x-2N_2\right)\ln|x|+{\rm O}(1)\quad \text{as}\quad |x|\to \infty, \\
 \ds r\pa_ru\to 2N_1-\frac{1}{2\pi}\int_{\mathbb{R}^2}\re^v\left(\re^u+1\right)\ud x,\quad \pa_\theta u\to 0\quad \text{as}\quad |x|\to \infty,\\
  \ds  r\pa_rv\to 2N_2-\frac{1}{2\pi}\int_{\mathbb{R}^2}\re^u\left(\re^v-1\right)\ud x,\quad \pa_\theta v\to 0\quad \text{as}\quad |x|\to \infty,
    \ea\right.\label{ss12}
 \eer
with $(r, \theta)$ the polar coordinates in $\mathbb{R}^2$.  More precisely, to establish \eqref{ss12}, we can argue exactly as for the proof of Theorem 1.1 in \cite{chta1}
by using Lemma 1.1 therein.

Therefore, by the integrability conditions,  and by recalling \eqref{ill},  we find that,
\ber\left\{\ba{lll}
 \ds\frac1\pi\Phi_1=\frac{1}{2\pi}\int_{\mathbb{R}^2}\re^v\left(\re^u+1\right)\ud x>2(N_1+1),\\[2mm]
 \ds\frac1\pi\Phi_2=\frac{1}{2\pi}\int_{\mathbb{R}^2}\re^u\left(\re^v-1\right)\ud x>2(N_2+1).
   \ea\right. \label{ss13}
\eer

Furthermore, we can use the above decay estimates at infinity to derive the following Pohozaev  identity:
\begin{lemma}
Every solution $(u, v)$ of \eqref{ss11} satisfies:
\ber
 &&\left(\frac{1}{2\pi}\int_{\mathbb{R}^2}\re^v\left(\re^u+1\right)\ud x-2(N_1+1)\right)\left(\frac{1}{2\pi}\int_{\mathbb{R}^2}\re^u\left(\re^v-1\right)\ud x-2(N_2+1)\right)+\frac{1}{\pi}\int_{\mathbb{R}^2}\re^{u+v}\ud x
 \nm\\
 &&=4(N_1+1)(N_2+1).\label{ss14}
\eer
\end{lemma}

{\it Proof.} By recalling  the identity:
\ber
 -\big[(x\cdot\nabla v)\Delta u+(x\cdot\nabla u)\Delta v\big]={\rm div} \big[x(\nabla u\cdot\nabla v)\big]-{\rm div} \big[(x\cdot\nabla v)\nabla u+(x\cdot\nabla u)\nabla v\big],
\eer
and by multiplying the first equation of  \eqref{ss11} by $x\cdot\nabla v$, and the second equation of \eqref{ss11} by  $x\cdot\nabla u$,  we can use  integration by parts
over  $\Omega\equiv  \left\{x\in \mathbb{R}^2\big|\, \delta<|x|<R\right\}$ for $0<\delta<R$ and obtain:
\ber
  \int_{\pa\Omega}(\nabla u\cdot\nabla v)x\cdot \nu \ud S-\int_{\pa\Omega}\big[(x\cdot\nabla v)\nabla u\cdot \nu +(x\cdot\nabla u)\nabla v\cdot \nu\big] \ud S\nm\\
  = \int_{\pa\Omega}\big[\re^v\left(\re^u+1\right)+\re^u\left(\re^v-1\right)\big]x\cdot \nu \ud S-\int_{\pa\Omega}\re^{u+v}x\cdot \nu \ud S-2\int_\Omega\left(\re^{u+v}+\re^v-\re^u\right)\ud x.\label{ss17'}
 \eer
Since after integration by parts, we also find that:
 \ber
  \ds -\int_{\partial\Omega}\nabla u\cdot\nu\ud S=\int_\Omega\re^v(\re^u+1)\ud x,\label{ss15'}\\
   \ds -\int_{\partial\Omega}\nabla v\cdot\nu\ud S=\int_\Omega\re^u(\re^v-1)\ud x,\label{ss16'}
 \eer
by combining the identities above, we obtain:
 \ber
   -\int_{\pa\Omega}\big[2(\nabla u\cdot \nu+\nabla v\cdot \nu)-(\nabla u\cdot\nabla v)x\cdot \nu+(x\cdot\nabla v)\nabla u\cdot \nu+(x\cdot\nabla u)\nabla v\cdot \nu\big] \ud S\nm\\
   = \int_{\pa\Omega}\big[\re^v\left(\re^u+1\right)+\re^u\left(\re^v-1\right)\big]x\cdot \nu \ud S-\int_{\pa\Omega}\re^{u+v}x\cdot \nu \ud S+2\int_\Omega\re^{u+v}\ud x.\label{ss18}
 \eer

Therefore, we can use the asymptotic  estimates \eqref{ss12}, and the logarithmal behavior of $u, v$ at the origin, to see that,  by letting $\delta\to 0$ and $R\to \infty$ in \eqref{ss18}, we can   obtain the desired identity \eqref{ss14}.
\hfill$\square$\\[1mm]

At this point, by straightforward calculation,  we can rewrite \eqref{ss14} in the following equivalent forms:
\ber
 &&\left(\frac{1}{2\pi}\int_{\mathbb{R}^2}\re^u\left(\re^v-1\right)\ud x-2(N_1+N_2+1)\right)\frac{1}{2\pi}\int_{\mathbb{R}^2}\re^{u+v}\ud x
 +\frac{N_1+1}{\pi}\int_{\mathbb{R}^2}\re^{u}\ud x\nm\\
 &&+\left(\frac{1}{2\pi}\int_{\mathbb{R}^2}\re^u\left(\re^v-1\right)\ud x-2(N_2+1)\right)\frac{1}{2\pi}\int_{\mathbb{R}^2}\re^{v}\ud x=0\label{ss15}
\eer
and
\ber
 &&\left(\frac{1}{2\pi}\int_{\mathbb{R}^2}\re^v\left(\re^u+1\right)\ud x-2(N_1+N_2+1)\right)\frac{1}{2\pi}\int_{\mathbb{R}^2}\re^{u+v}\ud x
-\frac{N_2+1}{\pi}\int_{\mathbb{R}^2}\re^{v}\ud x\nm\\
 &&-\left(\frac{1}{2\pi}\int_{\mathbb{R}^2}\re^v\left(\re^u+1\right)\ud x-2(N_1+1)\right)\frac{1}{2\pi}\int_{\mathbb{R}^2}\re^{u}\ud x=0.\label{ss16}
\eer

Consequently,  in view of \eqref{ss13},  we derive that  necessarily there must hold:
\ber
 2(N_2+1)<\frac1\pi\Phi_2=\frac{1}{2\pi}\int_{\mathbb{R}^2}\re^u\left(\re^v-1\right)\ud x<2(N_1+N_2+1)<\frac{1}{2\pi}\int_{\mathbb{R}^2}\re^v\left(\re^u+1\right)\ud x=\frac1\pi\Phi_1\label{ss17}
\eer
which express the necessary conditions claimed in \eqref{fl1} for the solvability of \eqref{ss11}.

Our next  goal in the following sections  is to show that, each of the condition above,  concerning  either $\Phi_1$ or $\Phi_2$  is also sufficient  for the solvability of \eqref{ss11}.

 \section{A perturbation approach and proof of Theorem \ref{th0}}\label{sec4}
\setcounter{equation}{0}\setcounter{lemma}{0}\setcounter{theorem}{0}
In this section we establish an existence theorem for the system \eqref{ss11} by means of a perturbation approach, from which  we can prove Theorem \ref{th0}.

For any $\vep>0$, we consider  the scaled functions:
 \ber
  u_\vep(x)=u(\vep x)+\ln\vep, \quad  v_\vep(x)=v(\vep x)+\ln\vep  \label{e3}
 \eer
satisfying:
  \ber\left\{\ba{lll}
  \ds -\Delta u_\vep=\re^{u_\vep+v_\vep}+\vep \re^{v_\vep}-4\pi N_1\delta_0\quad \text{in}\quad \mathbb{R}^2,\\[2mm]
   \ds -\Delta v_\vep=\re^{u_\vep+v_\vep}-\vep \re^{u_\vep}-4\pi N_2\delta_0\quad \text{in}\quad \mathbb{R}^2,\\[2mm]
   \ds    \re^{u_\vep},\,\re^{v_\vep},\,\re^{u_\vep+v_\vep} \in L^1(\mathbb{R}^2)
   \ea\right. \label{e4}
  \eer
  and
  \ber
 \ds \int_{\mathbb{R}^2}\re^{u_\vep+v_\vep}\ud x=\int_{\mathbb{R}^2}\re^{u+v}\ud x,\, \int_{\mathbb{R}^2}\re^{u_\vep}\ud x=\vep\int_{\mathbb{R}^2}\re^{u}\ud x,\,
 \ds \int_{\mathbb{R}^2}\re^{v_\vep}\ud x=\vep\int_{\mathbb{R}^2}\re^{v}\ud x.
  \eer

 Formally, by letting $\vep\to0^+$,  \eqref{e4}  can be viewed  as  a ``perturbation''  of the following problem:
    \ber\left\{\ba{lll}
  \ds -\Delta u_0=\re^{u_0+v_0}-4\pi N_1\delta_0\quad \text{in}\quad \mathbb{R}^2, \\[2mm]
   \ds -\Delta v_0=\re^{u_0+v_0}-4\pi N_2\delta_0\quad \text{in}\quad \mathbb{R}^2,\\[2mm]
   \ds    \re^{u_0},\,\re^{v_0},\,\re^{u_0+v_0} \in L^1(\mathbb{R}^2).
   \ea\right.\label{e5}
  \eer
To  classify  solutions of  the system  \eqref{e5},  simply we rewrite  it equivalently as follows:
 \ber\left\{\ba{lll}
 \ds  -\Delta(u_0+v_0)=2\re^{u_0+v_0}-4\pi(N_1+N_2)\delta_0\quad \text{in}\quad \mathbb{R}^2,\\[2mm]
 \ds  -\Delta(u_0-v_0)=-4\pi(N_1-N_2)\delta_0\quad \text{in}\quad \mathbb{R}^2,\\[2mm]
   \ds    \re^{u_0},\,\re^{v_0},\,\re^{u_0+v_0} \in L^1(\mathbb{R}^2),
  \ea\right. \label{e6}
 \eer
 and  so, by  the result   of  \cite{prta}  about  Liouville equation,
 we can explicitly  express    solutions of \eqref{e6} in terms of the free parameters $a\in \mathbb{C}$, $\lm>0$ and $c\in \mathbb{R}$, and
 (in complex notations) as follows:
 \berr\left\{\ba{lll}
\ds (u_0+v_0)(z)&=&\ln\frac{4(N_1+N_2+1)^2\lm|z|^{2(N_1+N_2)}}{(1+\lm|z^{N_1+N_2+1}+a|^2)^2},\\[2mm]
\ds  (u_0-v_0)(z)&=&(N_1-N_2)\ln|z|^2+c.
 \ea\right.
\eerr
We focus around the  following  (normalized) solution of \eqref{e5}:
  \ber \left\{\ba{lll}
 \ds u_0&=&\ln\frac{2(N_1+N_2+1)|z|^{2N_1}}{(1+|z|^{2(N_1+N_2+1)})},\\[3mm]
\ds  v_0&=&\ln\frac{2(N_1+N_2+1)|z|^{2N_2}}{(1+|z|^{2(N_1+N_2+1)})},
 \ea\right.\label{e7}
\eer
 and set,
  \ber
  \rho(r)&\equiv&2\re^{u_0+v_0}=\frac{8(N_1+N_2+1)^2r^{2(N_1+N_2)}}{(1+r^{2(N_1+N_2+1)})^2}, \quad r=|z|.\label{e8}
  \eer

We aim to  seek  solutions of \eqref{e4} in the form:
  \ber
  u_\vep(r)=u_0(r)+\vep u_1(r), \quad v_\vep(r)=v_0(r)+\vep v_1(r),  \label{e9}
  \eer
where  $(u_1, v_1)$ satisfies:
   \ber\left\{\ba{lll}
    \ds -\Delta u_1&=&\re^{u_0+v_0}\left(\frac{\re^{\vep(u_1+v_1)}-1}{\vep}\right)+\re^{v_0+\vep v_1},\\[2mm]
    \ds -\Delta v_1&=&\re^{u_0+v_0}\left(\frac{\re^{\vep(u_1+v_1)}-1}{\vep}\right)-\re^{u_0+\vep u_1}.
    \ea\right.\label{e10}
   \eer

Formally, by setting  $\vep\to0^+$, we see that  \eqref{e10}  yields to the following system:
      \ber\left\{\ba{lll}
    \ds -\Delta u_2&=&\re^{u_0+v_0}(u_2+v_2)+\re^{v_0},\\[2mm]
    \ds -\Delta v_2&=&\re^{u_0+v_0}(u_2+v_2)-\re^{u_0},
    \ea\right.\label{e11}
   \eer
or   equivalently:
     \ber\left\{\ba{lll}
    \ds -\Delta (u_2+v_2)&=&2\re^{u_0+v_0}(u_2+v_2)-\re^{u_0}+\re^{v_0}, \\[2mm]
   \ds  -\Delta  (u_2-v_2)&=&\re^{u_0}+\re^{v_0},
    \ea\right.\label{e12}
   \eer
with $u_0$ and $v_0$  specified in \eqref{e7}.

  We  are going to exhibit in the following   an explicit  radial  solution  $\big(u_2(r), v_2(r)\big)$ of \eqref{e12},  see  \eqref{e29}--\eqref{e30} below.

Therefore,  we use the decomposition:
 \ber
  u_1(r)=u_2(r)+u_3(r), \quad v_1(r)=v_2(r)+v_3(r)\label{e13}
 \eer
 with $u_3(r)$ and $ v_3(r)$ the error terms, going to zero in suitable norm as $\vep\to 0^+$.  So via \eqref{e9}, \eqref{e13} we determine a solution of the original problem
 \eqref{e10}, as soon as we solve:
  \ber\left\{\ba{lll} \label{e14}
  \ds \Delta u_3+\re^{u_0+v_0}\left(\frac{\re^{\vep(u_2+v_2+u_3+v_3)}-1}{\vep}-(u_2+v_2)\right)+\re^{v_0}(\re^{\vep(v_2+v_3)}-1)=0, \\[3mm]
  \ds \Delta v_3+\re^{u_0+v_0}\left(\frac{\re^{\vep(u_2+v_2+u_3+v_3)}-1}{\vep}-(u_2+v_2)\right)-\re^{u_0}(\re^{\vep(u_2+u_3)}-1)=0,
   \ea\right. \eer
   with the required integrability conditions.

Set
\ber\left\{\ba{lll}
  \ds P_1(u_3, v_3, \vep)&\equiv&\Delta u_3+\re^{u_0+v_0}\left(\frac{\re^{\vep(u_2+v_2+u_3+v_3)}-1}{\vep}-(u_2+v_2)\right)\\[2mm]
  &&+\re^{v_0}(\re^{\vep(v_2+v_3)}-1), \\[2mm]
  \ds  P_2(u_3, v_3, \vep)&\equiv&\Delta v_3+\re^{u_0+v_0}\left(\frac{\re^{\vep(u_2+v_2+u_3+v_3)}-1}{\vep}-(u_2+v_2)\right)\\[2mm]
   &&-\re^{u_0}(\re^{\vep(u_2+u_3)}-1),
   \ea\right.\label{e14'}
  \eer
and
 \ber
  P_i(0,0,0)\equiv0, \quad i=1, 2.  \label{e15}
 \eer

  Following  \cite{chaim}, it is natural to consider the operator $P_i$ to act between  the following spaces:
  \ber
  &&Y_\alpha \equiv \left\{w\in W^{2, 2}_{loc}(\mathbb{R}^2): \quad \Delta w\in X_\alpha, \quad  \frac{w}{(1+|z|^{1+\frac{\alpha}{2}})}\in L^2(\mathbb{R}^2)\right\}, \quad 0<\alpha<1,\label{e17}\\
  &&X_\alpha \equiv \Big\{w\in L^2_{loc}(\mathbb{R}^2): \quad (1+|z|^{2+\alpha})w^2\in L^1(\mathbb{R}^2)\Big\}, \quad 0<\alpha<1. \label{e16}
  \eer
Clearly, $Y_\alpha $ and $X_\alpha $ define  Hilbert spaces  equipped respectively  with  the following  scalar products:
 \ber
  &&(u, v)_{Y_\alpha}\equiv (\Delta u, \Delta v)_{X_\alpha}+\int_{\mathbb{R}^2} \frac{uv}{(1+|z|^{2+\alpha})}\ud z,\label{e19}\\
    &&(u, v)_{X_\alpha}\equiv\int_{\mathbb{R}^2}(1+|z|^{2+\alpha})uv\ud z\label{e18}
 \eer
 and induced  norms  denoted    respectively  by $\|\cdot\|_{Y_\alpha}$ and $\|\cdot\|_{X_\alpha}$.

 Moreover,  we let  $Y_{\alpha}^r$ and $X_{\alpha}^r$ respectively to  be the subspaces of $Y_\alpha$ and $X_{\alpha}$ consisting of   radial functions.

We recall the  following properties of functions in $Y_\alpha$,  whose proof can be found   in \cite{tabk}.
 \begin{lemma} {(\cite{tabk})} \label{lem1}
 Let $\alpha\in(0, 1)$, $w\in Y_\alpha$.

    (i) If  $w\in Y_\alpha$ is harmonic, then $w$ is a constant. \\

    (ii) The following estimates hold,\,\, $\forall\, w\in Y_\alpha$:
      \ber
      && |w(z)|\le C\|w\|_{Y_\alpha}\ln(1+|z|),\,\, \forall\, z\in \mathbb{R}^2, \label{e20'}  \\
       && \|\nabla w\|_{L^p}\le C_p\|w\|_{Y_\alpha}, \,\, \forall\, p>2 \label{e20''}
      \eer
      with $C$ and $C_p$ positive constants depending only  on $\alpha$ and $(\alpha, p)$, respectively.
 \end{lemma}
 \hfill $\square$

 It is easy to check that the   operator:
  \ber
  P:  Y_\alpha^r\times Y_\alpha^r\times\mathbb{R}\to X_\alpha^r\times X_\alpha^r, \label{e20}
    \eer
    with
   \ber
   P(u_3, v_3, \vep)\equiv\left(P_1(u_3, v_3, \vep),\,  P_2(u_3, v_3, \vep)\right) \label{e21}
   \eer
 is well defined, and it is continuous and differentiable  in $(0,0,0)$. Furthermore, we can  compute  the linearized operator $A$ at $(0, 0, 0)$,  namely:
 \ber
 A\equiv\frac{\pa P}{\pa (u_3, v_3)}(0, 0, 0):\quad   Y_\alpha^r\times Y_\alpha^r \to X_\alpha^r\times X_\alpha^r, \label{e22}
 \eer
and show that it  takes  the following form:
 \ber
A\begin{pmatrix}
 \phi\\\psi
\end{pmatrix}=
\begin{pmatrix}
 \Delta\phi+\re^{u_0+v_0}(\phi+\psi)\\
  \Delta\psi+\re^{u_0+v_0}(\phi+\psi)
\end{pmatrix},\quad \forall \, \phi,\psi \in Y_\alpha^r. \label{e23}
 \eer

To identify  the kernel of the operator $A$, we just need to solve  the problem:
\ber\left\{\ba{lll}
\ds \Delta(\phi+\psi)+\rho(\phi+\psi)=0,\\[2mm]
\ds \Delta(\phi-\psi)=0,\\[2mm]
 \ds \phi,\,\psi\in Y_\alpha^r
\ea\right. \label{e25}
\eer
with $\rho$   defined by \eqref{e8}.

To this purpose, we use  Lemma 3.4.20 and Corollary 3.4.21 in \cite{tabk},  for
  the operator  $L^r: Y^r_\alpha\to X_\alpha^r$  defined as follows:
 \ber
 L^rw\equiv \frac{\ud^2}{\ud r^2}w+\frac1rw+\rho w.\label{e26}
 \eer
More precisely, for
\ber
 \phi_0(r)\equiv\frac{1-r^{2(N_1+N_2+1)}}{1+r^{2(N_1+N_2+1)}}, \label{e27}
\eer
we have:
\begin{lemma}(\cite{tabk}) \label{lem2}
Let $\alpha\in(0, 1)$  and $N_1, N_2>0$.  We have

(i) $w\in Y_\alpha^r$ satisfies $L^r w=0$ if and only if $w \in \text{span}\{\phi_0\}$.

(ii) $L^r: Y^r_\alpha\to X_\alpha^r$ is onto. More precisely, for $f\in X_\alpha^r$, let
 \ber
  w(r)&=&\left(\phi_0(r)\ln r+\frac{2}{(N_1+N_2+1)(1+r^{2(N_1+N_2+1)})}\right)\int_0^r\phi_0(s)f(s)s\ud s\nn\\
  &&-\phi_0(r)\int_0^r\left(\phi_0(s)\ln s+\frac{2}{(N_1+N_2+1)(1+s^{2(N_1+N_2+1)})}\right)f(s)s\ud s,\label{e28}
 \eer
then $w\in Y_\alpha^r$ and $L^r w=f$. Moreover,
  \ber
   w(r)=-c_f\ln r+{\rm O}(1), \quad  w'(r)=-\frac{c_f}{r}+{\rm O}(1), \quad \text{as}\quad r\to \infty, \label{e28a}
  \eer
with
   \ber
    c_f\equiv\int_0^\infty\phi_0(t)f(t)t\ud t. \label{e28b}
   \eer
\end{lemma}
\hfill$\square$

 As a consequence, from  \eqref{e25} and  part $(i)$ of  Lemma \ref{lem1} and Lemma \ref{lem2},   we derive:
 \begin{corollary}
 \ber
  {\rm Ker} A= \left\{(w+c,\, w-c)\Big|\quad w\in {\rm span}\{\phi_0\}, \, c\in \mathbb{R}\right\}.\label{e31}
 \eer
  \end{corollary}
\hfill$\square$

Furthermore,
\begin{lemma}\label{lem3}
 The linearized operator $A$: $Y^r_\alpha\times Y^r_\alpha \to X^r_\alpha\times X^r_\alpha$ is onto.
\end{lemma}
 {\it Proof.} For any $(f_1, f_2)\in X^r_\alpha\times X^r_\alpha$,  we need to show that the system:
 \ber\left\{\ba{lll}
   \Delta \phi+\re^{u_0+v_0}(\phi+\psi)=f_1,\\[3mm]
   \Delta \psi+\re^{u_0+v_0}(\phi+\psi)=f_2
   \ea\right. \label{e32}\eer
 admits a solution in $Y^r_\alpha\times Y^r_\alpha$.

 To this purpose, we write  system \eqref{e32}  equivalently as  follows:
  \ber\left\{\ba{lll}
   \Delta (\phi+\psi)+\rho(\phi+\psi)=f_1+f_2,\\[3mm]
   \Delta (\phi-\psi)=f_1-f_2.
   \ea\right. \label{e33}\eer
 So,   by Lemma \ref{lem2},  we can solve  the first equation of  \eqref{e33} in $Y_\alpha^r$,  and obtain that,
 \berr
 &&(\phi+\psi)(r)\nn\\
 &&=\left(\phi_0(r)\ln r+\frac{2}{(N_1+N_2+1)(1+r^{2(N_1+N_2+1)})}\right)\int_0^r\phi_0(s)(f_1(s)+f_2(s))s\ud s\nn\\
  &&-\phi_0(r)\int_0^r\left(\phi_0(s)\ln s+\frac{2}{(N_1+N_2+1)(1+s^{2(N_1+N_2+1)})}\right)(f_1(s)+f_2(s))s\ud s.\label{e34}
 \eerr

While, as a particular solution of the second equation of  \eqref{e33}, we take:
\berr
(\phi-\psi)(r)=(\ln r)\int_0^rt(f_1(t)-f_2(t))\ud t-\int_0^rt(\ln t)(f_1(t)-f_2(t))\ud t
\eerr
and we easily check that $(\phi-\psi)(r) $  is well defined and moreover   $(\phi-\psi)(r)\in Y_\alpha^r$.

 At this point, we readily  derive  $(\phi, \psi)\in Y_\alpha^r\times Y_\alpha^r$  which  provides a solution for \eqref{e33}. \hfill$\square$\\[3mm]

By the expression of $u_0, v_0$ in \eqref{e7},  we see that $\re^{u_0}$ and $\re^{v_0}$ belong to $X_\alpha^r$.  Therefore,  by means of  the arguments  above,  we obtain
   a radial solution of   \eqref{e12}  as follows:
 \ber
 u_2(r)&=&\frac12\left(\phi_0(r)\ln r+\frac{2}{(N_1+N_2+1)(1+r^{2(N_1+N_2+1)})}\right)\int_0^r\phi_0(t)\left(\re^{u_0(t)}-\re^{v_0(t)}\right)t\ud t\nn\\
  &&-\frac12\phi_0(r)\int_0^r\left(\phi_0(t)\ln t+\frac{2}{(N_1+N_2+1)(1+s^{2(N_1+N_2+1)})}\right)\left(\re^{u_0(t)}-\re^{v_0(t)}\right)t\ud t\nn\\
  &&-\frac12(\ln r)\int_0^rt\left(\re^{u_0(t)}+\re^{v_0(t)}\right)t\ud t+\frac12\int_0^rt(\ln t)\left(\re^{u_0(t)}+\re^{v_0(t)}\right)\ud t,\label{e29}\\[3mm]
 v_2(r)&=&\frac12\left(\phi_0(r)\ln r+\frac{2}{(N_1+N_2+1)(1+r^{2(N_1+N_2+1)})}\right)\int_0^r\phi_0(t)\left(\re^{u_0(t)}-\re^{v_0(t)}\right)t\ud t\nn\\
  &&-\frac12\phi_0(r)\int_0^r\left(\phi_0(t)\ln t+\frac{2}{(N_1+N_2+1)(1+s^{2(N_1+N_2+1)})}\right)\left(\re^{u_0(t)}-\re^{v_0(t)}\right)t\ud t\nn\\
  &&+\frac12(\ln r)\int_0^rt\left(\re^{u_0(t)}+\re^{v_0(t)}\right)\ud t-\frac12\int_0^rt(\ln t)\left(\re^{u_0(t)}+\re^{v_0(t)}\right)\ud t.\label{e30}
 \eer

 Furthermore, the  values of the  finite  integrals:
 \ber
  &&\sigma_1\equiv\int_0^\infty\phi_0(t)\left(\re^{u_0(t)}-\re^{v_0(t)}\right) t\ud t,\label{e35}\\
  &&\sigma_2\equiv\int_0^\infty\left(\re^{u_0(t)}+\re^{v_0(t)}\right) t\ud t\label{e36}
 \eer
 characterize  the asymptotic behavior for $(u_2, v_2)$ as follows:
 \ber
  u_2(r)=-\frac{\sigma_1+\sigma_2}{2}\ln r+{\rm O}(1), \quad v_2(r)=-\frac{\sigma_1-\sigma_2}{2}\ln r+{\rm O}(1), \quad \text{as}\quad r\to \infty.\label{e37}
 \eer
\\[1mm]

Let $U_\alpha$ be the subspace of $Y_\alpha^r\times Y_\alpha^r$ orthogonal to ${\rm Ker} A$. Namely:
 \ber
 U_\alpha\equiv  ({\rm Ker} A)^\perp\subset Y_\alpha^r\times Y_\alpha^r,\label{e38}
 \eer
we establish the following existence result:
\begin{theorem}\label{thim}
 For fixed  $\alpha\in(0, 1)$, there exist  $\vep_0>0$ sufficiently small and   smooth functions $(u_{3,\vep}, v_{3, \vep}):\,(-\vep_0, \vep_0)\rightarrow U_\alpha$
 with $(u_{3,0}, v_{3, 0})=(0, 0)$ such that, by setting:
  \ber
 u^\vep(r)&\equiv&u_0\left(\frac r\vep \right)+\vep u_2\left(\frac r\vep \right)+\vep u_{3,\vep}\left(\frac r\vep\right)+\ln\frac1\vep, \label{e39}\\
 v^\vep(r)&\equiv& v_0\left(\frac r\vep \right)+\vep v_2\left(\frac r\vep \right)+\vep v_{3,\vep}\left(\frac r\vep\right)+\ln\frac1\vep, \label{e40}
 \eer
then  $(u^\vep(r),  v^\vep(r))$   defines   a radial  solution of the system \eqref{ss11}  with $(u_2, v_2)$ given by \eqref{e29}--\eqref{e30} and satisfying the asymptotic behavior \eqref{e37}.
 Moreover, as $\vep\to 0^+$, $\|u_{3,\vep}\|_{Y_\alpha^r}+\|v_{3,\vep}\|_{Y_\alpha^r}\to 0$ and  the following estimates hold:
  \ber
   &&|u_{3, \vep}(r)|+|v_{3, \vep}(r)|= \left(\|u_{3,\vep}\|_{Y_\alpha^r}+\|v_{3,\vep}\|_{Y_\alpha^r}\right)\ln(1+r),\quad \forall\, r>0,\label{e41}\\
   &&|u_{3, \vep}'(r)|+ |v_{3, \vep}'(r)|=\frac{1}{1+r}\left(\|u_{3,\vep}\|_{Y_\alpha^r}+\|v_{3,\vep}\|_{Y_\alpha^r}\right),\quad \forall\, r>0, \label{e42}\\
   &&\int_{\mathbb{R}^2}\re^{v^\vep}\left(\re^{u^\vep}+1\right)\ud x=4\pi(N_1+N_2+1)+{\rm o}(1), \label{e43}\\
   &&\int_{\mathbb{R}^2}\re^{u^\vep}\left(\re^{v^\vep}-1\right)\ud x=4\pi(N_1+N_2+1)+{\rm o}(1). \label{e44}
  \eer

 \end{theorem}
{\it Proof.}
Since  the linearized operator $A$ defined by \eqref{e22}  gives an isomorphism  between $U_\alpha$ and $X_\alpha^r\times X_\alpha^r$,
 we can use  the implicit function theorem (see e.g. \cite{nlbg}) for the operator $P: U_\alpha\times \mathbb{R}\to X_\alpha^r\times X_\alpha^r$
 around the point $(0,0,0)$. Therefore,  for sufficiently small $\vep_0>0$,  there exists
 a continuous function: $\vep\to \mathbf{\Psi}_\vep\equiv(u_{3, \vep}, v_{3, \vep})$ form  $ (-\vep_0, \vep_0)$ into a neighborhood of the origin in $U_\alpha$ such that
 $P(u_{3, \vep}, v_{3, \vep}, \vep)=0, \,\, \forall\,\vep\in(-\vep_0, \vep_0)$ and $(u_{3, 0}, v_{3, 0})=(0, 0)$.  As a consequence,  $(u^\vep(r),  v^\vep(r))$ defined by  \eqref{e39}--\eqref{e40} defines   a radial  solution
 of  the system \eqref{ss11}.

In particular, $\|u_{3,\vep}\|_{Y_\alpha^r}+\|v_{3,\vep}\|_{Y_\alpha^r}\to 0$ as $\vep\to 0^+$, and  we can use   \eqref{e20'}  together  with the equations \eqref{e14},  to deduce  the estimates \eqref{e41} and \eqref{e42}.

In view of \eqref{e39}--\eqref{e40}, we have:
\ber
 \int_{\mathbb{R}^2}\re^{v^\vep}\left(\re^{u^\vep}+1\right)\ud x =\int_{\mathbb{R}^2}\re^{v_0(r)+\vep v_2(r)+\vep v_{3,\vep}(r)}\left(\re^{u_0(r)+\vep u_2(r)+\vep u_{3,\vep}(r)}+\vep\right)\ud x, \label{e45}
\eer
and by the asymptotic estimates for $u_2, v_2, u_{3, \vep},$ and $v_{3, \vep}$ established above, we  conclude that,
\ber
\lim\limits_{\vep\to0}\int_{\mathbb{R}^2}\re^{v^\vep}\left(\re^{u^\vep}+1\right)\ud x
&=&\int_{\mathbb{R}^2}\re^{u_0(r)+v_0(r)}\ud x\nn\\
&=&8\pi(N_1+N_2+1)^2\int_0^\infty \frac{r^{2(N_1+N_2)}}{\left(1+r^{2(N_1+N_2+1)}\right)^2}r\ud r\nm\\
&=&4\pi(N_1+N_2+1), \label{e46}
\eer
as claimed.

Similarly, we obtain:
\ber
\lim\limits_{\vep\to0}\int_{\mathbb{R}^2}\re^{u^\vep}\left(\re^{v^\vep}-1\right)\ud x=4\pi(N_1+N_2+1),  \label{e47}
\eer
and the proof is completed. \hfill $\square$\\[1mm]

{\bf The proof of Theorem \ref{th0}}.

At this point, by means of  Theorem \ref{thim} we can establish easily Theorem \ref{th0}.

In fact,   by using  \eqref{i48}--\eqref{i49} (considered with the lower sign), together with the relations \eqref{re1}--\eqref{re3},
we  find that,
 \ber
  &&f^{(1),\vep}_{12}=\frac{8\pi^2}{k^2}|\phi_2^\vep|^2\left(|\phi_1^\vep|^2+c^2\right)=2\sigma^2\re^{v^\vep}\left(\re^{u^\vep}+1\right),\label{4.50}\\
  &&f^{(2),\vep}_{12}=\frac{8\pi^2}{k^2}|\phi_1^\vep|^2\left(|\phi_2^\vep|^2-c^2\right)=2\sigma^2\re^{u^\vep}\left(\re^{v^\vep}-1\right).\label{4.51}
 \eer

 On the other hand, if we use the decomposition  \eqref{e39} and \eqref{e40} for $u^\vep$ and $v^\vep$  together with \eqref{e7} and the given decay estimates for the error term, we easily check that,
 as $\vep\to0^+$,
 \ber
  &&\re^{u^\vep+v^\vep} \rightharpoonup 4\pi (N_1+N_2+1)\delta_0\quad \text{weakly in the sense of measure in }\,\, \mathbb{R}^2,\label{4.52}\\
  &&\re^{u^\vep},\,\, \re^{v^\vep}\to 0\quad \text{in}\quad L^1(\mathbb{R}^2).\label{4.53}
 \eer

Therefore, by taking into account the normalization \eqref{2.78a}  and \eqref{i66}, from \eqref{4.50}--\eqref{4.51},  we readily  deduce   \eqref{2.54} and \eqref{2.55}.
Clearly, also  \eqref{2.56} and  \eqref{2.57} are a direct consequence of \eqref{e43}--\eqref{e44}.  Finally,  \eqref{2.53a} follows from \eqref{e7} and
  the given decay estimate of the error terms.

 Theorem \ref{th0}  is established.

\section {The radial problem}\label{sec5}
\setcounter{equation}{0}\setcounter{lemma}{0}\setcounter{theorem}{0}\setcounter{corollary}{0}

In this section we collect some general properties about radial solutions of  \eqref{ss11}.

 More precisely, for $u=u(r), v=v(r)$, we  consider  the following problem:
 \ber\left\{\ba{lll}\label{e1'}
\ds-(ru')'=r\re^v\left(\re^u+1\right), \, r>0, \\[2mm]
\ds-(rv')'=r\re^u\left(\re^v-1\right),\, r>0, \\[2mm]
\ds u(r)-2N_1\ln r= {\rm O}(1) \quad \text{as}\quad r\to 0^+, \\[2mm]
\ds v(r)-2N_2\ln r= {\rm O}(1)\quad \text{as}\quad r\to 0^+.
   \ea\right.\eer

We start to point out the following extended version of  Pohozaev identity \eqref{ss14} valid   in the radial setting:
 \begin{proposition} \label{pro1}
 For any  solution $(u, v)$ of  the  problem  \eqref{e1'},  we have:
 \ber
 &&2\int_0^rt\re^{u(t)+v(t)}\ud t\nm\\
 &&=4(N_1+1)(N_2+1)-\big(ru'(r)+2\big)\big(rv'(r)+2\big)-r^2\left(\re^{u(r)+v(r)}+\re^{v(r)}-\re^{u(r)}\right),\label{e4'}\\
 &&2\int_0^rt\re^{u(t)}\ud t=4N_1(N_2+1)-ru'(r)(rv'(r)+2)+ r^2\left(\re^{u(r)+v(r)}+\re^{v(r)}-\re^{u(r)}\right),\label{e5'}\\
 &&2\int_0^rt\re^{v(t)}\ud t=rv'(r)(ru'(r)+2)-4N_2(N_1+1)+ r^2\left(\re^{u(r)+v(r)}+\re^{v(r)}-\re^{u(r)}\right).\label{e6'}
 \eer
 \end{proposition}
{\it Proof.}
  By integration and by using \eqref{e1'} we find:
   \ber
    ru'(r)=2N_1-\int_0^rt\re^{v(t)}\left(\re^{u(t)}+1\right)\ud t,\label{e7'}\\
    rv'(r)=2N_2-\int_0^rt\re^{u(t)}\left(\re^{v(t)}-1\right)\ud t. \label{e8'}
   \eer
Multiplying the first equation in \eqref{e1'} by $rv'$, and the second equation in   \eqref{e1'} by $ru'$,  and  after integrating by parts, we obtain:
\ber
 &&4N_1N_2-r^2u'(r)v'(r)\nm\\
 &&= r^2\left(\re^{u(r)+v(r)}+\re^{v(r)}-\re^{u(r)}\right)-2\int_0^rt\re^{u(t)+v(t)}\ud t-2\int_0^rt\left(\re^{v(t)}-\re^{u(t)}\right)\ud t.\label{e9'}
\eer

Combining \eqref{e7'}--\eqref{e9'},  we obtain  \eqref{e4'}--\eqref{e6'}.  \hfill $\square$\\[1mm]

Interestingly, for radial solutions of \eqref{e1'},  the following integrability properties hold:

\begin{lemma}\label{lem5.1}
Let $(u,v)$ be a  solution of \eqref{e1'}, then we have:
\ber\int_0^\infty r\re^{u+v}\ud r<\infty, \quad \int_0^\infty r\re^u\ud r<\infty. \label{s12}
\eer
\end{lemma}

{\it Proof. }We argue by contradiction and suppose $\int_0^rt\re^{u+v}\ud t\to \infty$ as $r\to \infty$.
 We see that $ru'(r)$ is monotone decreasing and diverges to $-\infty$ as $r\to\infty$. As a consequence, $\frac{u(r)}{\ln r} \to -\infty$ as $r\to\infty$,
and so,  $\int_0^\infty r\re^u\ud r<\infty$. At this point we can use \eqref{e8'} to conclude that also, $\lim\limits_{r\to \infty}\frac{v(r)}{\ln r}=-\infty$, and so:
 $\int_0^\infty r\re^{u+v}\ud r<\infty$, a  contradiction.

Similarly, we get that, $\int_0^\infty r\re^{u}\ud r<\infty$.  Indeed, if by  contradiction, we suppose that,  $\int_0^rt\re^u\ud t\to \infty$ as $r\to \infty$, then  from  \eqref{e8'},  we find  that $rv'\to \infty $ as $r\to  \infty$. As a consequence,
 $\int_0^rt\re^v\ud t \to \infty$ as $r\to\infty$ and  in view of \eqref{e7'},  we have  $ru'\to -\infty$ as $r\to\infty$. In other words, $\frac{u(r)}{\ln r}\to-\infty$ and so
   $\int_0^\infty r\re^u\ud r<\infty$,  a contradiction.   \hfill$ \square$
\\[2mm]

In view of Lemma \ref{lem5.1}, to obtain a (radial) solution of \eqref{ss11}, we only need to
ensure the integrability of $\re^v$. To this purpose, let
 \ber u(r)=U(r)+2N_1\ln r, \quad  v(r)=V(r)+2N_2\ln r. \label{s0}
 \eer
From now on  we use  the following notation:
 \ber\left\{\ba{lll}\label{s00}
 \ds F_1(r)\equiv\int_0^rt^{2N_2+1}\re^{V(t)}\left(t^{2N_1}\re^{U(t)}+1\right)\ud t=\int_0^rt\re^{v(t)}\left(\re^{u(t)}+1\right)\ud t, \\[3mm]
  \ds F_2(r)\equiv\int_0^rt^{2N_1+1}\re^{U(t)}\left(t^{2N_2}\re^{V(t)}-1\right)\ud t=\int_0^rt\re^{u(t)}\left(\re^{v(t)}-1\right)\ud t;
   \ea\right.\eer
and  we  consider the initial value problem for the (regular) functions $(U(r), V(r))$ associated to \eqref{e1'}. Namely:
 \ber\left\{\ba{lll}
  \ds-(rU')'=r^{2N_2+1}\re^V\left(r^{2N_1}\re^U+1\right),  \\[2mm]
  \ds-(rV')'=r^{2N_1+1}\re^U\left(r^{2N_2}\re^V-1\right),  \\[2mm]
  \ds U(0)=\alpha_1,\, V(0)=\alpha_2, \\[2mm]
  \ds  U'(0)=V'(0)=0,\label{s1}
   \ea\right.\eer
with $(\alpha_1, \alpha_2)\in \mathbb{R}^2$.

For the  initial value problem \eqref{s1},  the following holds:
\begin{proposition} \label{prop1}
For any  $(\alpha_1, \alpha_2)\in \mathbb{R}^2$, the  initial value problem \eqref{s1} admits a unique  global solution,
 which depends continuously on the initial data.
\end{proposition}
{\it Proof.}  By standard ODE  techniques, we see that the  system \eqref{s1} admits a unique local solution $(U(r),V(r))$ defined  in the interval $(0, r_0)$ for some $r_0>0$.
By integration,  for every  $r\in(0, r_0)$ we have:
 \ber
  rU'(r)=-F_1(r)<0  \quad \text{and} \quad rV'(t)=-F_2(r). \label{s3}
 \eer
In particular, we see that  $U(r)$ is strictly  decreasing in $(0, r_0)$, and so  $U(r)\le U(0)=\alpha_1$,   for $0<r<r_0$.
While, by \eqref{s3}, for $0<r<r_0$ we get,
\ber
 V(r)=V(0)-\int_0^r\frac1tF_2(t)\ud t\le V(0)+\int_0^r\frac1t\int_0^ts^{2N_1+1}\re^{U(s)}\ud s\ud t
  \le V(0)+\frac{\re^{U(0)}r_0^{2(N_1+1)}}{4(N_1+1)^2}. \label{s4}
\eer

By  \eqref{s3},   we see that,   for every  $0<r<r_0$ the following holds:
 \ber
  |U'(r)|\le \left|\frac1rF_1(r)\right|\le C \quad \text{and}\quad  |V'(r)|\le \left|\frac1rF_2(r)\right|\le C.\label{s6}
 \eer
 for some constant $C>0$ depending only on the initial data $(\alpha_1, \alpha_2)$.

Hence,  for all $0\le r< r_0$, we have:
\ber
 |U(r)| \le |U(0)|+Cr_0 \quad  \text{and}\quad  |V(r)|\le |V(0)|+Cr_0. \label{s8}
\eer

 At this point, we can use a  standard unique continuation argument to conclude that the initial value problem \eqref{s1}
 admits a unique  global solution, which depends continuously on the initial data $(\alpha_1, \alpha_2)\in \mathbb{R}^2$. \hfill $\square$

\begin{remark} \label{rmk51}
In view of Proposition \ref{prop1}, we know that every solution of the radial problem \eqref{e1'} corresponds (via \eqref{s0}) to  a solution of
the initial value problem \eqref{s1} and viceversa, in other words,  problem \eqref{e1'}  and  \eqref{s1} are equivalent.

Furthermore, from  \eqref{s3}  we see that,  in case $F_1(\infty)<\infty$, then the following  estimate holds:
 \[U(r)>-F_1(\infty)\ln r-C, \,\,\forall \, r\ge 1, \]
with $C>0$  a suitable constant, and it  implies  that necessarily $F_1(\infty)>2(N_1+1)$ (as already pointed out in Section \ref{sec3}).
 With this information, we can show  in turn that a similar estimate holds for $U+V$, namely:
  \ber
  U(r)+V(r)\ge -\big(F_1(\infty)+F_2(\infty)\big)\ln r-C, \quad \forall\, r\ge 1,
  \eer
and by Lemma \ref{lem5.1}, we conclude that:  $F_1(\infty)+F_2(\infty)>2(N_1+N_2+1)$.  As a consequence,
  \ber
   r^{2(N_1+N_2+1)}\re^{U(r)+V(r)}\to 0, \quad       r^{2(N_1+1)}\re^{U(r)}\to 0    \quad\text{as}\,\, r\to \infty.   \label{s8a}
  \eer
\end{remark}

 To proceed further, we point out some qualitative information about solutions of \eqref{e1'}.  First  of all, from \eqref{s0} and \eqref{s8a}, we see that
 $u(r)\to-\infty$ as $r\to0^+$ and $u(r)+2\ln r\to-\infty$ as $r\to\infty$, so  both $u(r)$ and $u(r)+2\ln r$ admit a unique maximum attained respectively
 at the (unique) values: $0<t_1<\tilde{t}_1$ satisfying: $F_1(t_1)=2N_1$ and  $F_1(\tilde{t}_1)=2(N_1+1)$, see \eqref{e7'}.
 Furthermore,   we can also conclude that,   in case the maximum value of $u(r)$ (or $u(r)+2\ln r $) is positive then  $u(r)$ (or $u(r)+2\ln r $) must vanish  exactly twice.

 Less obvious is the behavior of $v(r)$,  and we can describe it  in terms of  the value of $F_2(\infty)$ as follows:
 \begin{proposition} \label{prop4.3}

Let $(u, v)$  satisfy  \eqref{e1'}. There  holds:

(i) $v(r)$ vanishes at least once at a value $r_0>0$ which corresponds to the global minimum point of $F_2(r)$. Furthermore, we have that, $v(r)<0, \,\,\forall\, r\in(0, r_0)$.

(ii) $F_2(r)$ can attain the value $2N_2$ at most once and:

a) If $F_2(\infty)<2N_2$, then $F_2(r)<2N_2,\,\,\forall\, r>0$,  $v$ is strictly increasing (to $\infty$) with unique zero at $r_0$, and so $F_2(r)$ is strictly decreasing in $(0, r_0)$ and strictly
increasing in $[r_0, \infty)$;

b) If $F_2(\infty)>2N_2$, then $v$ vanishes exactly twice at $r_0<t_0$, and  $r_0$ corresponds to the global minimum point of $F_2(r)$ while  $t_0$ corresponds to the global maximum  point of $F_2(r)$. Moreover, $F_2$ is decreasing in
$[0, r_0]\cup[t_0, \infty)$ and increasing in  $(r_0, t_0)$.  Furthermore, $v$ admits a unique maximum at $r_*\in (r_0, t_0)$ satisfying $F_2(r_*)=2N_2$.

 \end{proposition}

{\it Proof.}
To  prove (i), we first observe that $v(r)<0$ for $r$ sufficiently close to $r=0$. So   the set $\Lambda\equiv\{r>0| \, v(t)\le 0, \,\, \forall\, t\in(0, r]\}$ is not empty.

{\bf Claim 1}: $\Lambda$ is bounded from above and $r_0=\sup \Lambda$.

To establish Claim 1, we argue by contradiction and assume that $v(r)\le 0$, $\forall\, r>0$. As a consequence, $F_2(r)<0, \,\, \forall\, r>0$ and so the function
$V(r)=v(r)-2N_2\ln r$ is strictly increasing and must satisfy: $r^{2N_2}\re^{V(r)}=\re^{v(r)}\le 1,\, \forall\, r>0$. In particular, $r^{2N_2}\re^{V(0)}\le 1$ for all $r>0$,  and this  is clearly  impossible.

At this point, if we set $r_0=\sup \Lambda$,  we  find that: $v(r_0)=0$ and $F_2(r_0)<0$. So, for $\delta >0$ sufficiently small we still have $F_2(r)<0, \forall\, r\in[0, r_0+\delta]$.
Thus, $v$ is increasing in $(0, r_0+\delta)$ and in particular $v(t)<0, \,\,\forall\, t\in (0, r_0)$,  while $v(r)>0$ as $r\in (r_0, r_0+\delta)$. By those information
we conclude  the proof of Claim 1 and also see that $r_0$ defines a strictly local  minimum for $F_2(r)$. We shall show below that
 actually $r_0$ is a global minimum for $F_2(r)$.

But first we  establish (ii).  We see that,  $F_2(r)<0$ for $r\in(0, r_0]$ and so we consider the set
 \[\Lambda_1=\left\{ r>0\big|\, \forall t\in (0, r):\,  F_2(t)<2N_2\right\}\supset (0, r_0]. \]
In case $\Lambda_1$ is unbounded above, we see that in this case $F_2(r)$ never attains the value $2N_2$.  Hence we assume that, $\tau_1=\sup\Lambda_1<\infty$.
Consequently, $F_2(\tau_1)=2N_2$ and $F_2(r)<2N_2,\,\,  \forall\, r\in [0, \tau_1)$. Furthermore, we know that $\tau_1>r_0$ and $v$ is strictly increasing in $[0, \tau_1]$.  As a consequence, $v(\tau_1)>0(=v(r_0))$.

{\bf Claim 2}:  $F_2(r)>2N_2, \,\,  \forall\, r>\tau_1$.

To establish Claim 2, we observe first that, for $\vep>0$ sufficiently small, we have: $v(r)>0,\,\, \forall\, r\in [\tau_1, \tau_1+\vep]$, and so $F_2(r)>2N_2,\,\, \forall\, r\in(\tau_1, \tau_1+\vep]$.
We argue again by contradiction and assume that,
\[\tau_2=\sup \left\{r>\tau_1, \,\forall\, t\in(\tau_1, r)\big|\, F_2(r)>2N_2\right\}<\infty.\]
Hence, $\tau_2>\tau_1$, $F_2(\tau_2)=2N_2$ and $F_2(r)>2N_2, \,\, \forall\, r\in(\tau_1, \tau_2)$.   In particular, $v(r)$ is strictly decreasing in $(\tau_1, \tau_2)$
and $v(\tau_2)<0$. In other words, we have obtained the following:
 \ber
 && 0<\tau_1<\tau_2: F_2(\tau_1)=F_2(\tau_2)=2N_2 \,\,\text{and}\,\,v(\tau_2)<0<v(\tau_1),\label{31*}\\
 && v(\tau_2)<v(\tau)<v(\tau_1),\,\, \forall\, r\in (\tau_1, \tau_2).\label{32*}
  \eer
We show next  that \eqref{31*}--\eqref{32*} are  impossible (i.e. $v$ cannot ``oscillate''). Indeed, by using Pohozaev identity \eqref{e4'}  together with \eqref{e7'} and \eqref{e8'} respectively at $r=\tau_1$ and $r=\tau_2$,  we have:
\berr
\left[2(N_1+1)-F_1(\tau_1)\right]\left[2(N_2+1)-F_2(\tau_1)\right]+2\int_0^{\tau_1}t\re^{u(t)+v(t)}\ud t +\tau_1^2\left(\re^{u(\tau_1)+v(\tau_1)}+\re^{v(\tau_1)}-\re^{u(\tau_1)}\right)\nm\\
=\left[2(N_1+1)-F_1(\tau_2)\right]\left[2(N_2+1)-F_2(\tau_2)\right]+2\int_0^{\tau_2}t\re^{u(t)+v(t)}\ud t +\tau_2^2\left(\re^{u(\tau_2)+v(\tau_2)}+\re^{v(\tau_2)}-\re^{u(\tau_2)}\right),
\eerr
and since $F_2(\tau_1)=F_2(\tau_2)=2N_2$, we obtain:
\berr
-2\int_0^{\tau_1}t\re^{v(t)}\ud t +\tau_1^2\left(\re^{u(\tau_1)+v(\tau_1)}+\re^{v(\tau_1)}-\re^{u(\tau_1)}\right)
=-2\int_0^{\tau_2}t\re^{v(t)}\ud t +\tau_2^2\left(\re^{u(\tau_2)+v(\tau_2)}+\re^{v(\tau_2)}-\re^{u(\tau_2)}\right).
\eerr
Consequently,  since   $v(\tau_2)<0<v(\tau_1)$,  we derive:
\ber
2\int_{\tau_1}^{\tau_2}t\re^{v(t)}\ud t&=&\tau_2^2\re^{u(\tau_2)}\left(\re^{v(\tau_2)}-1\right)+\tau_2^2\re^{v(\tau_2)}-\tau_1^2\re^{u(\tau_1)}\left(\re^{v(\tau_1)}-1\right)-\tau_1^2\re^{v(\tau_1)}\nm\\
&<& \tau_2^2\re^{v(\tau_2)}-\tau_1^2\re^{v(\tau_1)}.
\eer

 On the other hand,   $v$ is strictly decreasing in $[\tau_1, \tau_2]$, and so we have:
\[
 (\tau_2^2-\tau_1^2)\re^{v(\tau_2)}<2\int_{\tau_1}^{\tau_2}t\re^{v(t)}\ud t\le \tau_2^2\re^{v(\tau_2)}-\tau_1^2\re^{v(\tau_1)},
 \]
which  implies that:   $\tau_1^2(\re^{v(\tau_1)}-\re^{v(\tau_2)})< 0$, a contradiction to the  fact that:  $v(\tau_1)>v(\tau_2)$.  Therefore, Claim 2 is established, as well as (ii).

If $F_2(\infty)<2N_2$, by virtue of (ii)  we see that necessarily:  $F_2(r)<2N_2, \,\,  \forall\, r>0$. Hence $v$ is strictly increasing and $v(r)\to \infty$ as $r\to \infty$, and so,
$r_0$ is the only zero of $v$. So,  $F_2(r)$ is strictly decreasing in $(0, r_0)$ and strictly increasing in $(r_0, \infty)$, and consequently $r_0$ corresponds to its global minimum point.
On the other hand, if $F_2(\infty)>2N_2$, then by (ii), we know that there exists a unique value $r_*>0$ such that, $F_2(r)<2N_2,\,\,\forall\, r\in(0, r_*)$ and $F_2(r)>2N_2$
for $r>r_*$ (see Claim 2) and $r_*>r_0$ (see Claim 1).  Consequently, $v$ is strictly increasing in $(0, r_*)$ and strictly decreasing in $(r_*, \infty)$ and $v$
attains its maximum value at $r_*$. Furthermore:  $\frac{v(r)}{\ln r}\to \big(2N_2-F_2(\infty)\big)<0$ as $r\to \infty$, so  we see that $v(r)\to -\infty$ as $r\to \infty$,  and so $v$ must vanish once at $t_0\in (r_*, \infty)$.
This implies that $F_2(r)$ is strictly decreasing in $(0, r_0]\cup[t_0, \infty)$, while it is increasing in $(r_0, t_0)$, so $r_0$ and $t_0$ corresponds respectively to the  global minimum and
global maximum points for $F_2(r)$.

The case $F_2(\infty)=2N_2$ enters as a  border line case,  and for $v(r)$ and $F_2(r)$ both behaviors a) and b) could occur. However,  in any case, $r_0$ would always correspond to a global minimum point
for $F_2(r)$, and the proof is completed. \hfill $\square$\\[1mm]

From Proposition \ref{prop4.3}  we can derive easily some useful information about $(U,V)$  solutions  of \eqref{s1}. Clearly, $U$ is strictly decreasing (to $-\infty$), namely: $U(r)\le U(0),\,\, \forall\, r>0$.
Concerning $V$, we have:
\begin{corollary}\label{crl3}
 Let $(U, V)$ be a solution of \eqref{s1}, we have:

 (i) If $F_2(\infty)\le0$, then $F_2(r)<0,\,\, \forall\, r>0$ and $V$  is  strictly increasing (to $\infty$);

 (ii) If $F_2(\infty)>0$, then there exists a unique $s_0>r_0$ such that $F_2(s_0)=0$,  and $F_2(s)<0, \,\, \forall \,s\in(0, s_0)$, while  $F_2(s)>0$ for $s>s_0$. In particular,
 $s_0$ is the (unique) maximum point of $V$.
\end{corollary}

Our next goal is to provide the following a priori estimates.

\begin{proposition}\label{prp4}
For a solution $(U, V)$ of \eqref{s1},  the following estimates hold:
\ber
 &&i)\quad \frac{1}{2(N_1+1)}\left(\max\limits_{r\in [0,\infty)}r^{2(N_1+1)}\re^{U(r)}\right)\le \int_0^\infty t^{2N_1+1}\re^{U(t)}\ud t \nm\\
 &&\qquad \qquad \qquad \qquad \qquad \qquad \qquad \qquad \, \le 2\left[2(N_1+1)+\re^{U(0)-\frac{N_1+1}{N_2}V(0)}\right],\quad \label{n1}\\
&&ii) \quad F_2(r)<2(N_1+N_2+1), \label{n1'}\\
&&\quad \big(F_1(r)-2(N_1+N_2+1)\big) \int_0^rt^{2(N_1+N_2)+1}\re^{U(t)+V(t)}\ud t+r^{2(N_1+N_2+1)}\re^{U(r)+V(r)}>0,\qquad    \label{n1''}
\eer
$\forall r\ge0$;   and in particular,
\ber
&&\int_0^\infty t^{2(N_1+N_2)+1}\re^{U(t)+V(t)}\ud t \le 2\left[3(N_1+1)+N_2+\re^{U(0)-\frac{N_1+1}{N_2}V(0)}\right],\label{n2}\\
&&\max\limits_{r\in [0,\infty)}\left(r^{2(N_1+N_2+1)}\re^{U(r)+V(r)}\right)\le \left[3(N_1+1)+N_2+\re^{U(0)-\frac{N_1+1}{N_2}V(0)}\right]^2.\label{n2'}
\eer

\end{proposition}

{\it Proof.}  Recall that, $U(r)\le U(0)$ and $V(r)\ge V(0),\,\,\forall r\in [0, r_0]$,  with $r_0>0$ as given by
 part (i) of Proposition \ref{prop4.3}.  Recall  also  that,
 \ber
  F_2(r_0)<0\quad \text{and}\quad  F_2(r)\ge F_2(r_0), \,\, \forall\, r\ge0. \label{n3'}
 \eer
 Therefore, by means of such information and Jensen's inequality, we have:
 \berr
 \int_0^{r_0}t^{2N_1+1}\re^{U(t)}\ud t&\ge&\int_0^{r_0}t^{2(N_1+N_2)+1}\re^{U(t)+V(t)}\ud t\nm\\
 &=& \int_0^{r_0}t^{2(N_1+N_2)+1}\re^{\frac{2(N_1+N_2)+1}{2N_1+1}U(t)}\re^{V(t)-\frac{2N_2}{2N_1+1}U(t)}\ud t\nm\\
 &\ge& r_0\re^{V(0)-\frac{2N_2}{2N_1+1}U(0)}\dashint_0^{r_0}\left(t^{2N_1+1}\re^{U(t)}\right)^{\frac{2(N_1+N_2)+1}{2N_1+1}}\ud t\nm\\
 &\ge& r_0\re^{V(0)-\frac{2N_2}{2N_1+1}U(0)}\left(\dashint_0^{r_0}t^{2N_1+1}\re^{U(t)}\ud t\right)^{\frac{2(N_1+N_2)+1}{2N_1+1}},
 \eerr
which implies
\ber
 \int_0^{r_0}t^{2N_1+1}\re^{U(t)}\ud t\le r_0\re^{U(0)-\frac{2N_1+1}{2N_2}V(0)}. \label{n3a}
\eer
 By recalling that, $0=v(r_0)=V(r_0)+2N_2\ln r_0>V(0)+2N_2\ln r_0$, we see that $r_0<\re^{-\frac{V(0)}{2N_2}}$ and from \eqref{n3a} we obtain the following:
\ber
\int_0^{r_0}t^{2N_1+1}\re^{U(t)}\ud t\le  \re^{U(0)-\frac{N_1+1}{N_2}V(0)}.\label{n3b}
\eer
Clearly, from \eqref{n3b} we find:
 \ber
  |F_2(r_0)|\le \int_0^{r_0}t^{2N_1+1}\re^{U(t)}\ud t\le \re^{U(0)-\frac{N_1+1}{N_2}V(0)}.\label{n3c}
 \eer

To  prove \eqref{n1},   we  use   the first equation in \eqref{s1} and integration by parts,  to find:
 \ber
  &&2(N_1+1)\int_0^rt^{2N_1+1}\re^{U(t)}\ud t\\
  &&=r^{2(N_1+1)}\re^{U(r)}+\int_0^rt^{2N_1+1}\re^{U(t)}F_1(t)\ud t\nm\\
  &&=r^{2(N_1+1)}\re^{U(r)}+\int_0^rt^{2N_1+1}\re^{U(t)}\left(\int_0^ts^{2(N_1+N_2)+1}\re^{U(s)+V(s)}\ud s-\int_0^ts^{2N_1+1}\re^{U(s)}\ud s\right)\nm\\
  &&\quad +\frac12\left(\int_0^rt^{2N_1+1}\re^{U(t)}\ud t\right)^2+\int_0^rt^{2(N_1+N_2)+1}\re^{U(t)+V(t)}\int_0^t s\re^{v(s)}\ud s\ud t\nm\\
  &&\ge r^{2(N_1+1)}\re^{U(r)}+\int_0^rt^{2N_1+1}\re^{U(t)}F_2(t)\ud t+\frac12\left(\int_0^rt^{2N_1+1}\re^{U(t)}\ud t\right)^2\nm\\
    &&\ge  r^{2(N_1+1)}\re^{U(r)}+ F_2(r_0)\int_0^rt^{2N_1+1}\re^{U(t)}\ud t+\frac12\left(\int_0^rt^{2N_1+1}\re^{U(t)}\ud t\right)^2,\label{n3d}
 \eer
 which implies:
 \ber
\frac{r^{2(N_1+1)}\re^{U(r)}}{2(N_1+1)}\le \int_0^rt^{2N_1+1}\re^{U(t)}\ud t\le   2\big[2(N_1+1)+|F_2(r_0)|\big], \quad \forall\, r>0.\label{n4}
 \eer

 Thus,  we can use \eqref{n3c} into \eqref{n4}, and by letting  $r\to \infty$, we obtain \eqref{n1}.

 Similarly, we can prove \eqref{n2}. Indeed, by simple integration by parts,   we have:
 \berr
 &&2(N_1+N_2+1)\int_0^rt^{2(N_1+N_2)+1}\re^{U(t)+V(t)}\ud t\nm\\
 &&=r^{2(N_1+N_2+1)}\re^{U(r)+V(r)}+\int_0^rt^{2(N_1+N_2)+1}\re^{U(t)+V(t)}(F_1(t)+F_2(t))\ud t\nm\\
 &&=r^{2(N_1+N_2+1)}\re^{U(r)+V(r)}+\left(\int_0^rt^{2(N_1+N_2)+1}\re^{U(t)+V(t)}\ud t\right)^2\nm\\
 &&\quad +\int_0^rt^{2(N_1+N_2)+1}\re^{U(t)+V(t)}\left(\int_0^t(s^{2N_2+1}\re^{V(s)}-s^{2N_1+1}\re^{U(s)})\ud s\right)\ud t.\label{s39a}
 \eerr
 As a consequence,
 \ber
  &&2(N_1+N_2+1)\int_0^rt^{2(N_1+N_2)+1}\re^{U(t)+V(t)}\ud t\nm\\
  &&\ge r^{2(N_1+N_2+1)}\re^{U(r)+V(r)}+\left(\int_0^rt^{2(N_1+N_2)+1}\re^{U(t)+V(t)}\ud t\right)^2 \nm\\
  &&\quad -\int_0^rt^{2(N_1+N_2)+1}\re^{U(t)+V(t)}\ud t\int_0^rt^{2N_1+1}\re^{U(t)}\ud t,\label{s39}
\eer
and
\ber\label{s39b}
&&2(N_1+N_2+1)\int_0^rt^{2(N_1+N_2)+1}\re^{U(t)+V(t)}\ud t < r^{2(N_1+N_2+1)}\re^{U(r)+V(r)}\nm\\
&&\quad +\left(\int_0^rt^{2(N_1+N_2)+1}\re^{U(t)+V(t)}\ud t\right)^2
+\int_0^rt^{2(N_1+N_2)+1}\re^{U(t)+V(t)}\int_0^rt^{2N_2+1}\re^{V(t)}\ud t.
\eer

From \eqref{s39b} we readily derive the  inequality claimed in \eqref{n1''}.  On the other hand, by setting
\ber
A(r)\equiv\int_0^rt^{2(N_1+N_2)+1}\re^{U(t)+V(t)}\ud t, \label{s40}
\eer
from \eqref{s39} we get the following  inequality:
\ber
 &&\left[A(r)-\left(N_1+N_2+1+\frac12\int_0^rs^{2N_1+1}\re^{U(s)}\ud s\right)\right]^2\nm\\
 &&\le \left(N_1+N_2+1+\frac12\int_0^rs^{2N_1+1}\re^{U(s)}\ud s\right)^2-r^{2(N_1+N_2+1)}\re^{U(r)+V(r)}.\label{s41}
\eer

Consequently, from   \eqref{s41},  we obtain, $\forall \, r>0$:
\ber
&& r^{2(N_1+N_2+1)}\re^{U(r)+V(r)}\le \left(N_1+N_2+1+\frac12\int_0^rs^{2N_1+1}\re^{U(s)}\ud s\right)^2\label{s42}
\eer
and
\ber
\int_0^r t^{2(N_1+N_2)+1}\re^{U(t)+V(t)}\ud t< 2\left(N_1+N_2+1+\frac12\int_0^r s^{2N_1+1}\re^{U(s)}\ud s\right).\label{s43}
\eer

At this point, from \eqref{s43} we readily derive \eqref{n1'}. Furthermore,  we may also conclude the desired estimates \eqref{n2}--\eqref{n2'},   by letting $r\to \infty$ and by using  \eqref{n1}. \hfill $\square$

\begin{remark} \label{rmk52}
Note in particular that from \eqref{n1'}  and \eqref{s8a} we obtain the necessary  condition:
\ber
F_2(\infty)<2(N_1+N_2+1)<F_1(\infty)\le \infty,\label{s43a}
\eer
 already derived for general finite-energy  solutions (i.e. not necessarily radial) in Section \ref{sec3}.
\end{remark}

 \section{Existence of radial solutions and the proof of Theorem \ref{th1} and Theorem \ref{th2}}\label{sec6}
 \setcounter{equation}{0} \setcounter{remark}{0}
  \setcounter{lemma}{0}\setcounter{corollary}{0}

 In this section we focus our attention to the radial solvability of \eqref{ss11} and  carry out the proof of  Theorem \ref{th1} and Theorem \ref{th2}.

 To this end,    by  Lemma \ref{lem5.1} and via \eqref{s0}, we need to identify the
 initial conditions for    the Cauchy problem \eqref{s1}  which yields to  a solution  satisfying:
 \ber
  \int_0^\infty t\re^{v(t)}\ud t=\int_0^\infty t^{2N_2+1}\re^{V(t)}\ud t<\infty. \label{z1}
 \eer

 According to the discussion in Section 3 (valid for solutions not necessarily radially symmetric),  we know that,
 \begin{lemma}\label{lemz1}
 A solution of \eqref{s1} satisfies \eqref{z1} if and only if
  \ber
   F_2(\infty)>2(N_2+1). \label{z2}
  \eer
 \end{lemma}

Actually, in the radial case,  the statement of Lemma \ref{lemz1}  can be easily checked.

To this purpose, we recall that, as $r\to \infty$,
  \ber
 \frac{V(r)}{\ln r}\to -F_2(\infty) \quad \text{and}\quad \frac{U(r)}{\ln r}\to -F_1(\infty)<-2(N_1+N_2+1), \label{6.2a}
  \eer
  (see Remark \ref{rmk52}). While it is clear that \eqref{z2} implies \eqref{z1}, on the contrary  if \eqref{z1} holds then we  can conclude only that:
 $F_2(\infty)\ge 2(N_2+1)$. In particular,  if \eqref{z1} holds, then,
    \ber
    \frac{U(r)+V(r)}{\ln r}\to -(F_1(\infty)+F_2(\infty))<-2(N_1+N_2+1)-2(N_2+1), \quad \text{as}\quad r\to \infty.\label{6.2b}
    \eer
To check that actually,  if $F_2(\infty)=2(N_2+1)$ then  \eqref{z1} fails,  we  use \eqref{6.2a} and \eqref{6.2b}, to see that
\ber
 V(r)+2(N_2+1)\ln r=V(0)+(F_2(\infty)-F_2(r))\ln r+\int_0^r(\ln t)F'_2(t)\ud t={\rm O}(1),\,\text{as}\quad  r\to \infty,\label{6.2c}
 \eer
and  so $V$ cannot satisfy \eqref{z1}.
\begin{remark}
Since for a solution $(U, V)$ of \eqref{s1} satifying \eqref{z1} we have that, $F_1(\infty)<\infty$ and \eqref{6.2b} holds, then we can use those facts, together with
\eqref{z2} and \eqref{6.2a} to show (as in \eqref{6.2c}) that,
\ber
 U(r)=-F_1(\infty)\ln r+{\rm O}(1), \quad V(r)=-F_2(\infty)\ln r+{\rm O}(1),\quad \text{as}\quad r\to \infty,\label{6.2d}
\eer
and (by recalling \eqref{e7'} and \eqref{e8'}):
\ber
U'(r)=-\frac{F_1(\infty)}{r} +{\rm o}(1), \quad V'(r)=-\frac{F_2(\infty)}{r} +{\rm o}(1),\quad \text{as}\quad r\to \infty. \label{6.2e}
\eer
\end{remark}

More importantly, we observe that not all initial data yield to solutions satisfying \eqref{z1} (or equivalently \eqref{z2}).
To see this fact, we use some new estimates specific of    solutions  of  \eqref{s1}   satisfying \eqref{z1},  which also
have the advantage  to depend only  on $U(0)$.

\begin{lemma}\label{lemz2}
Let $(U,V)$ be a solution of \eqref{s1} satisfying the integrability condition \eqref{z1},   we have:
\ber
 &&r^{2(N_2+1)}\re^{V(r)}\le 4(N_1+1)(N_2+1), \,\forall\,r>0,\label{s45}\\
 &&r^{2(N_1+1)}\re^{U(r)}\le C(N_1, N_2) \left(1+\re^{U(0)}\right),\forall\, r>0,\label{s46}\\
 &&\int_0^\infty r^{2(N_1+N_2)+1}\re^{U(r)+V(r)}\ud r\le C(N_1, N_2)\re^{\frac{U(0)}{N_1+1}}\left(1+\re^{U(0)}\right)^{\frac{1}{N_1+1}},\label{s46'}
\eer
with  $C(N_1, N_2)$ a suitable  positive constant depending only on $N_1, N_2$.
\end{lemma}
{\it Proof.}
In view of \eqref{z1}, we know that $v(r)+2\ln r=V(r)+2(N_2+1)\ln r\to -\infty$ as $r\to 0^+$ and as
$r\to\infty$. Therefore, $V(r)+2(N_2+1)\ln r$ attains its  maximum value, say at $r_2>0$, with $F_2(r_2)=2(N_2+1)$.
By the information provided by Proposition \ref{prop4.3}, we also know that such maximum point is
unique and $V(r_2)-2N_2\ln r_2=v(r_2)>0$.  Thus, by \eqref{e8'} and  the  Pohozaev identity \eqref{e4'} applied  at $r=r_2>0$ we find:
 \ber
  &&2\int_0^{r_2}r^{2(N_2+N_2)+1}\re^{U(r)+V(r)}\ud r+r_2^{2(N_1+1)}\re^{U(r_2)}\left(r_2^{2N_2}\re^{V(r_2)}-1\right)+r_2^{2(N_2+1)}\re^{V(r_2)}\nm\\
  &&=4(N_1+1)(N_2+1),
 \eer
and  we readily deduce that,
 \ber
  r_2^{2(N_2+1)}\re^{V(r_2)}\le 4(N_1+1)(N_2+1),
 \eer
and \eqref{s45} is established.

Next, let $r_1>0$ be the unique maximum point of $u(r)+2\ln r$.  Namely,  $U(r_1)+2(N_1+1)\ln r_1=\max\{U(r)+2(N_1+1)\ln r\}$ and so, $F_1(r_1)=2(N_1+1)$. Pohozaev's  identity \eqref{e4'}  at $r=r_1>0$ gives:
 \ber
&&2\int_0^{r_1}r^{2(N_1+N_2)+1}\re^{U(r)+V(r)}\ud r+r_1^{2(N_1+N_2+1)}\re^{U(r_1)+V(r_1)}+r_1^{2(N_2+1)}\re^{V(r_1)}-r_1^{2(N_1+1)}\re^{U(r_1)}\nm\\
 &&= 4(N_1+1)(N_2+1),
 \eer
or equivalently,
 \ber
 r_1^{2(N_1+1)}\re^{U(r_1)}&=&2\int_0^{r_1}r^{2(N_1+N_2)+1}\re^{U(r)+V(r)}\ud r+r_1^{2(N_2+1)}\re^{V(r_1)}-4(N_1+1)(N_2+1)\nm\\
 &&+r_1^{2(N_1+N_2+1)}\re^{U(r_1)+V(r_1)}.\label{s47}
 \eer
 Therefore,  if we use \eqref{s47} and the fact that,
 \ber
 \int_0^{r_1}r^{2(N_1+N_2)+1}\re^{U(r)+V(r)}\ud r\le F_1(r_1)=2(N_1+1),
 \eer
 we have:
 \ber
  r_1^{2(N_1+1)}\re^{U(r_1)}&\le& 2N_1+r_1^{2(N_1+N_2+1)}\re^{U(r_1)+V(r_1)}\nm\\
  &\le& 2N_1+4(N_1+1)(N_2+1)r_1^{2N_1}\re^{U(r_1)}. \label{s52}
 \eer
In case there holds: $r_1^2\le 8(N_1+1)(N_2+1)$  then,
 \ber
 r_1^{2(N_1+1)}\re^{U(r_1)}\le C(N_1,N_2)\re^{U(0)},\label{s53}
 \eer
 with $C(N_1,N_2)=\big(8(N_1+1)(N_2+1)\big)^{N_1+1}$.

On the other hand, if $r_1^2>8(N_1+1)(N_2+1)$,  then we can use  \eqref{s52}  to obtain that,
 $r_1^{2(N_1+1)}\re^{U(r_1)}\le 4N_1$,  and so in all cases we see that   \eqref{s46} holds.

 Next  we prove \eqref{s46'}. For any $R>0$, we use  \eqref{s45} and \eqref{s46} to find:
 \ber
 \int_0^\infty t^{2(N_1+N_2)+1}\re^{U(t)+V(t)}\ud t &=& \int_0^Rt^{2(N_1+N_2)+1}\re^{U(t)+V(t)}\ud t+\int_R^\infty t^{2(N_1+N_2)+1}\re^{U(t)+V(t)}\ud t\nm\\
 &\le& C(N_1, N_2)\left(\int_0^Rt^{2N_1-1}\re^{U(t)}\ud t+\int_R^\infty\frac{1+\re^{U(0)}}{t^3}\right)\nm\\
 &\le& C(N_1, N_2)\left(\re^{U(0)}R^{2N_1}+\frac{1+\re^{U(0)}}{R^2}\right). \label{s54a}
 \eer
At this point, by minimizing with respect to $R>0$  the right hand side of \eqref{s54a} (namely,
by taking $R=\left(\frac{1+\re^{U(0)}}{N_1\re^{U(0)}}\right)^{\frac{1}{2(N_1+1)}}$  in \eqref{s54a}),  we obtain:
 \ber
   \int_0^\infty t^{2(N_1+N_2)+1}\re^{U(t)+V(t)}\ud t \le C(N_1, N_2)\re^{\frac{U(0)}{N_1+1}}\left(1+\re^{U(0)}\right)^{\frac{1}{N_1+1}},
 \eer
and \eqref{s46'} is also established. \hfill $\square$

\begin{corollary}\label{corl62}
There exists $a_0\equiv a_0(N_1, N_2)$ such that if $(U,V)$ is the solution of \eqref{s1} with $U(0)=\alpha_1<a_0$, then
\eqref{z1}(or \eqref{z2}) fails.
\end{corollary}
{\it Proof.} Indeed, if by contradiction we suppose that $F_2(\infty)>2(N_2+1)$, then we would have:
$2(N_2+1)\le \int_0^\infty t^{2(N_1+N_2)+1}\re^{U(t)+V(t)}\ud t,$ and  by means of the estimates \eqref{s46'}
we could certainly violate such an inequality by letting $U(0)\to -\infty$.  \hfill $\square$\\[1mm]

Inspired by the shape of the solution constructed in Section \ref{sec4}, we prove the following:

\begin{theorem} \label{thfr}
For any given $L\in \mathbb{R}$ and $\gamma\in\big(2(N_2+1),\, 2(N_1+N_2+1)\big)$(or $E>0$), there exists $\alpha_\gamma=\alpha_\gamma(L)$  (or $\alpha=\alpha_E(L)$)
 such that, the (unique) solution $(U, V)$ of \eqref{s1} with  $(U(0), V(0))=\left(\alpha_\gamma,\frac{(2N_2+1)\alpha_\gamma-L}{2N_1+1}\right)$
 satisfies:
 \ber
 \gamma=F_2(\infty)=\int_0^\infty r^{2N_1+1}\re^{U(r)}\left(r^{2N_2+1}\re^{V(r)}-1\right)\ud r.
 \eer
 Similarly,    if we fix $\gamma\in\big(2(N_1+N_2+1), \infty\big)$, then an analogous statement hold with
  \ber
 \gamma=F_1(\infty)=\int_0^\infty r^{2N_2+1}\re^{V(r)}\left(r^{2N_1+1}\re^{U(r)}+1\right)\ud r.
 \eer
 While,  if  $(U(0), V(0))=\left(\alpha_E,\frac{(2N_2+1)\alpha_E-L}{2N_1+1}\right)$, then it holds:
 \ber
 N(N-1)\sigma^3k\int_0^\infty (t\re^{u(t)}+t\re^{v(t)})\ud t= N(N-1)\sigma^3k\int_0^\infty\left(t^{2N_1+1}\re^{U(t)}+t^{2N_2+1}\re^{V(t)}\right)=E. \label{s66s}
 \eer
Furthermore, by setting  $\beta_i\equiv\frac12F_i(\infty)-N_i$ ($i=1, 2$), then $\beta_1>N_2+1$ and $1<\beta_2<N_1+1$, and
  \ber
   0<(\beta_1-1)(\beta_2-1)&<&(N_1+1)(N_2+1),\label{s68}\\
    \beta_2(\beta_1-1)&>&N_2(N_1+1), \label{s69}\\
    0<\beta_1(\beta_2-1)&<&N_1(N_2+1). \label{s70}
  \eer
\end{theorem}
\begin{remark}\label{rmk62}
The role of the parameters $\beta_1$ and $\beta_2$ introduced in Theorem \ref{thfr} can be justified better in terms of the pair $(u(r), v(r))$:
$u(r)=U(r)+2N_1\ln r$ and $v(r)=V(r)+2N_2\ln r,$
 which defines a radial solution of \eqref{ss11} satisfying:
  \ber
   &&u(r)=-2\beta_1\ln r+{\rm O}(1),\quad   v(r)=-2\beta_2\ln r+{\rm O}(1), \quad \text{as}\quad r\to \infty, \label{de1}\\
    &&u'(r)=-\frac{2\beta_1}{r}+{\rm o}(1),\quad   v'(r)=-\frac{2\beta_2}{r}+{\rm o}(1),  \quad \text{as}\quad r\to \infty,\label{de2}
  \eer
  as we can easily derive from \eqref{6.2d} and \eqref{6.2e}.  In other words, $2\beta_1$ and $2\beta_2$ identify the power of decay at infinity respectively
   of $\re^u$  and $\re^v$ (and their derivatives), consistently with \eqref{ss12}.
\end{remark}

To establish Theorem \ref{thfr},  we use a blow-up argument together with the information provided by Corollary \ref{corl62}.

\begin{theorem} \label{thbp} For any given $L\in \mathbb{R}$,  let $(U_n(r), V_n(r))$ be a sequence of solutions for \eqref{s1} with initial data
 $(U_n(0), V_n(0))$ satisfying
 \ber
(2N_2+1)U_n(0)-(2N_1+1)V_n(0)\to L,  \quad \text{as}\quad   n\to \infty. \label{5.1}
 \eer
Let
\ber
F_{2,n}(r)\equiv \int_0^r t^{2(N_1+N_2)+1}\re^{U_n(t)+V_n(t)}\ud t-\int_0^r t^{2N_1+1}\re^{U_n(t)}\ud t,
\eer
the following holds:
 \ber
 \text{if} \quad U_n(0)\to\infty, \quad\text{then}\quad F_{2, n}(\infty)\to 2(N_1+N_2+1),\quad \text{as}\quad n\to \infty. \label{5.1a}
 \eer
 \end{theorem}

{\it  Proof.}  The proof relies on a  blow-up analysis.

In view of Corollary \ref{crl3} we let  $0<t_n\le \infty$, be  such that $F_{2,n}(t)<0$ for all $t\in(0, t_n)$.

 {\bf Claim.} $0<t_n\le 1$  and so:   $F_{2,n}(t_n)=0,$ $V_n(t_n)=\max V_n$ and $F_{2,n}(t)>0$ for $t>t_n$.

We argue by contradiction and suppose that  $F_{2,n}(r)<0, \, \forall r\in[0, 1]$.  Hence,  $V_n(r)$ is increasing in $[0, 1]$,
and  we have:
\ber
0<V_n(r)-V_n(0)=-F_{2,n}(r)\ln r+\int_0^r(\ln t)  F_{2,n}'(t)\ud t,\,\, \forall\, r\in(0,1]. \label{5.3}
\eer
In particular, from \eqref{5.3} we find that,  $\int_0^r \ln t F_{2,n}'(t)\ud t>0$, $\forall r\in[0, 1]$.

Thus, by recalling that $U_n(r)\le U_n(0)$, $V_n(r)$ is increasing in $[0, 1]$,  and by using  \eqref{5.1}, we have:
\ber
&&\int_0^r\left(\ln\frac1t\right)  t^{2N_1+1}\re^{U_n(t)}\ud t\nm\\
&&>\int_0^r}\left(\ln\frac1t\right) t^{2(N_1+N_2)+1}\re^{U_n(t)+V_n(t)\ud t\nm\\
&&=\int_0^r\left(\ln\frac1t\right)  t^{2(N_1+N_2)+1}\re^{\frac{2(N_1+N_2+1)}{2N_1+1}U_n(t)+V_n(t)-\frac{2N_2+1}{2N_1+1}U_n(t)}\ud t\nm\\
&&\ge\int_0^r\left(\ln\frac1t\right)  t^{2(N_1+N_2)+1}\re^{\frac{2(N_1+N_2+1)}{2N_1+1}U_n(t)}\re^{V_n(t)-\frac{2N_2+1}{2N_1+1}U_n(0)}\ud t\nm\\
&&\ge\re^{-\frac{L}{2N_1+1}}(1+{\rm o}(1))\int_0^r\left(\ln\frac1t\right)  t^{2(N_1+N_2)+1}\re^{\frac{2(N_1+N_2+1)}{2N_1+1}U_n(t)} \ud t,\label{5.5}
\eer
as $n\to \infty$.
On the other hand, by  H\"{o}lder's inequality,  we get:
\berr
&&\int_0^r\left(\ln\frac1t\right) t^{2N_1+1}\re^{U_n(t)}\ud t\nm\\
&&\le \left(\int_0^r\left(\ln\frac1t\right)t^{2(N_1+N_2)+1}\re^{\frac{2(N_1+N_2+1)}{2N_1+1}U_n(t)} \ud t\right)^{\frac{2N_1+1}{2(N_1+N_2+1)}}
\left(\int_0^r\left(\ln\frac1t\right) t^{\frac{2N_1+1}{2N_2+1}}\ud t\right)^{\frac{2N_2+1}{2(N_1+N_2+1)}},
\eerr
which we can use together with \eqref{5.5} to obtain:
\berr
\left(\int_0^r\left(\ln\frac1t\right)  t^{2N_1+1}\re^{U_n(t)}\ud t\right)^{\frac{2N_2+1}{2N_1+1}}\le \re^{\frac{L}{2N_1+1}}(1+{\rm o}(1))\left(\int_0^r\left(\ln\frac1t\right)  t^{\frac{2N_1+1}{2N_2+1}}\ud t\right)^{\frac{2N_2+1}{2N_1+1}},\, \forall\, r\in[0, 1],\quad
\eerr
as $n\to \infty$.
 As a consequence, $\forall r\in[0, 1]$ and $\forall n\in \mathbb{N}$, we conclude that:
\ber
\int_0^r\left(\ln\frac1t\right)  t^{2(N_1+N_2)+1}\re^{U_n(t)+V_n(t)}\ud t<\int_0^r\left(\ln\frac1t\right)  t^{2N_1+1}\re^{U_n(t)}\ud t\le C, \label{5.8}
\eer
with  $C>0$  a suitable  constant depending only on $L$, $N_1$ and $N_2$.
 Therefore, from  \eqref{5.3} and \eqref{5.8} we deduce also  that,
\berr
0<V_n(r)-V_n(0)\le C,\,\, \forall\, r\in [0, 1].
\eerr

 To obtain the desired contradiction,  we use a blow-up analysis for the scaled functions:
\ber
 \tilde{U}_n(r)=U_n(s_nr)-U_n(0), \quad  \tilde{V}_n(r)=V_n(s_nr)-V_n(0), \quad r\in\left[0,\frac{1}{s_n}\right],\label{5.10}
\eer
with
\ber
s_n=\re^{-\frac{U_n(0)+V_n(0)}{2(N_1+N_2+1)}}. \label{5.11}
\eer
By means of \eqref{5.1}, we see that, as $n\to \infty$,
\berr
1&=&s_n^{2(N_1+N_2+1)}\re^{U_n(0)+V_n(0)}\nm\\
&=&\re^{-\frac{L}{2N_1+1}}(1+{\rm o}(1))s_n^{2(N_1+N_2+1)}\re^{\frac{2(N_1+N_2+1)}{2N_1+1}U_n(0)}\nm\\
&=&\re^{-\frac{L}{2N_1+1}}(1+{\rm o}(1))(s_n^{2N_1+1}\re^{U_n(0)})^{\frac{2(N_1+N_2+1)}{2N_1+1}},
\eerr
which implies  that,
\ber
s_n^{2N_1+1}\re^{U_n(0)}={\rm O}(1)  \label{5.12}
\eer
and similarly  we get,
\ber
 s_n^{2N_2+1}\re^{V_n(0)}={\rm O}(1).  \label{5.13}
\eer

Furthermore, $(\tilde{U}_n, \tilde{V}_n)$  can be casted as radially symmetric solutions in $D_n\equiv \left\{x\in R^2\big| \, |x|<\frac{1}{s_n}\right\}$ of the system:
  \ber\left\{\ba{lll}\label{r17}
\ds-\Delta\tilde{U}_n=|x|^{2(N_1+N_2)}\re^{\tilde{U}_n+\tilde{V}_n}+\vep_{2,n}|x|^{2N_2}\re^{\tilde{V}_n},\quad x\in D_n\\[2mm]
\ds-\Delta\tilde{V}_n=|x|^{2(N_1+N_2)}\re^{\tilde{U}_n+\tilde{V}_n}-\vep_{1,n}|x|^{2N_1}\re^{\tilde{U}_n},\quad x\in D_n
 \ea\right.\eer
with
\ber
 \vep_{1, n}\equiv s_{n}^{2(N_1+1)}\re^{U_n(0)}\to 0, \quad \vep_{2, n}\equiv s_{n}^{2(N_2+1)}\re^{V_n(0)}\to 0,   \quad \text{as}\,\,n\to \infty \label{r19}
\eer
 (see \eqref{5.12}--\eqref{5.13}) and in $D_n$  the following holds:
 \ber
 \tilde{U}_n\le \tilde{U}_n(0)=0, \quad  \tilde{V}_n\le C \quad\text{and}\,\, \tilde{V}_n(0)=0.\label{r16}
\eer

In particular, we have:
\ber
-\Delta(\tilde{U}_n+\tilde{V}_n)=2|x|^{2(N_1+N_2)}\re^{\tilde{U}_n+\tilde{V}_n}+\vep_{2,n}|x|^{2N_2}\re^{\tilde{V}_n}-\vep_{1,n}|x|^{2N_1}\re^{\tilde{U}_n},
\eer
 with $\tilde{U}_n(0)+\tilde{V}_n(0)=0$ and $\tilde{U}_n+\tilde{V}_n\le C$ in $D_n$.

 Furthermore, by using \eqref{5.1} with $U_n(0)\to +\infty$ into the estimates \eqref{n2},  or by using  estimates  similar to those provided in \eqref{5.5} (simply by dropping the term $\ln\frac1t$),
 we find a suitable constant $C>0$:
\ber
\int_0^{\frac{1}{s_n}}r^{2(N_1+N_2)}\re^{\tilde{U}_n(r)+\tilde{V}_n(r)}\ud r=\int_0^1r^{2(N_1+N_2)}\re^{U_n(r)+V_n(r)}\ud r<C.
\eer

Therefore, we are in a position to use  standard elliptic estimates together with an Harnack type inequality (see e.g. Corollary 5.2.9 in \cite{tabk}, and corresponding application  to
Liouville type equations discussed therein)  and obtain that,  along a subsequence, the following holds:
\ber
 \tilde{U}_n+\tilde{V}_n\to \xi, \quad  \text{uniformly in}\quad C_{loc}^2(\mathbb{R}^2), \quad  \text{as}\quad n\to \infty\label{r24}
\eer
with $\xi$  a radial solution for the Liouville problem:
\ber\left\{\ba{lll} \label{r25}
\ds-\Delta\xi=2r^{2(N_1+N_2)}\re^{\xi} \quad \text{in}\quad \mathbb{R}^2,\\[2mm]
\ds\int_{\mathbb{R}^2} r^{2(N_1+N_2)+1}\re^{\xi(r)}\ud r \le C, \\[2mm]
  \ds\xi(0)=0.
\ea\right.\eer

By the classification result  in \cite{CL1,CL2} and \cite{prta}, we know the explicit expression for $\xi=\xi(r)$ and  in particular that  it  satisfies:
\ber
 \int_0^\infty r^{2(N_1+N_2)+1}\re^{\xi(r)} \ud r=2(N_1+N_2+1). \label{r27}
\eer
As a consequence, for any $\vep>0$, there exist  $R_\vep>0$ and $n_\vep>0$ such that,
\ber
 \int_0^{R_\vep}r^{2(N_1+N_2)+1}\re^{\tilde{U}_n(r)+\tilde{V}_n(r)}\ud r\ge 2(N_1+N_2+1)-\vep, \quad \forall \,n\ge n_\vep.
\eer

 On the other hand,  since  $F_{2, n}(t)\le 0$ for $0\le t\le 1$,  we see that,
 \ber
  \int_0^{R_\vep}r^{2(N_1+N_2)+1}\re^{\tilde{U}_n(r)+\tilde{V}_n(r)}\ud r&=&\int_0^{s_nR_\vep}r^{2(N_1+N_2)+1}\re^{U_n(r)+V_n(r)}\ud r\nm\\
  &<&\int_0^{s_nR_\vep}r^{2N_1+1}\re^{U_n(r)}\ud r\nm\\
  &\le&\frac{\re^{U_n(0)}}{2(N_1+1)}s_n^{2(N_1+1)}R_\vep^{2(N_1+1)}\to 0 \,\,\text{as}\,\, n\to\infty,
 \eer
which  leads to  the desired  contradiction, and the Claim is established.

 Therefore, we can use the estimates \eqref{5.3} for $r=t_n$, and as above deduce that,
 \ber
 0<V_n(t_n)-V_n(0)=\int_0^{t_n}(\ln t) F'_{2,n}(t)\ud t\le C.
 \eer

Hence,  for the scaled functions $(\tilde{U}_n, \tilde{V}_n)$  now  we can claim that,
 \ber
 \tilde{U}_n(r)\le \tilde{U}_n(0)=0,\quad \tilde{V}_n(r)\le V_n\left(\frac{t_n}{s_n}\right)\le C,  \,\tilde{V}_n(0)=0,\label{6.37}\\
 \int_0^{\infty} r^{2(N_1+N_2)+1}\re^{U_n+V_n}\le C,\quad \int_0^{\infty} r^{2N_1+1}\re^{U_n}\le C, \label{6.38}
 \eer
 where the estimates \eqref{6.38} hold  with a constant $C>0$ depending only on $N_1,N_2$ and $L$, and they can be derived by using \eqref{5.1} with $U_n(0)\to \infty$,
 into the estimates \eqref{n1} and \eqref{n2'},  as  follows:
\ber
&&\int_0^\infty t^{2N_1+1}\re^{U_n(t)}\ud t\le 4(N_1+1)+2\re^{U_n(0)-\frac{N_1+1}{N_2}V_n(0)}=4(N_1+1)+{\rm o}(1),\label{6.39}\\
  &&\int_0^\infty t^{2(N_1+N_2)+1}\re^{U_n(t)+V_n(t)}\ud t\le 6(N_1+1)+2N_2+2\re^{U_n(0)-\frac{N_1+1}{N_2}V_n(0)}\nm\\
  &&\qquad\qquad\qquad\qquad\qquad\qquad\quad= 6(N_1+1)+2N_2+{\rm o}(1),\label{6.40}
\eer
as  $n\to \infty$.

As above, we can  show that:
\ber
0<\frac{t_n}{s_n}<C,  \quad  \text{for some}\quad  C>0. \label{6.41}
\eer
Indeed, by  assuming  that (along a subsequence),  $\frac{t_n}{s_n}\to \infty$, then as before,  we can carry out a blow-up argument in $ [0, \frac{t_n}{s_n})$
to get a  contradiction.

 Therefore, in view of \eqref{6.37}, \eqref{6.38},  \eqref{6.41}  as above by  well-known elliptic and Harnack estimates, we  obtain  that,  along a subsequence,
\ber
 \tilde{U}_n\to U, \quad \tilde{V}_n\to V\quad \text{in}\quad C_{loc}^2(\mathbb{R}^2)\label{r24'}
\eer
 with $(U, V)$ a radial solution  of  the problem:
\ber\left\{\ba{lll}\label{r25'}
 \ds -\Delta U=|x|^{2(N_1+N_2)}\re^{U+V}\quad \text{in}\quad  \mathbb{R}^2,\\[2mm]
 \ds -\Delta V=|x|^{2(N_1+N_2)}\re^{U+V}\quad \text{in}\quad  \mathbb{R}^2,\\[2mm]
 \ds \int_{\mathbb{R}^2} |x|^{2(N_1+N_2)}\re^{U+V}\ud x <\infty,\quad  U(0)=0=V(0).
 \ea\right.
\eer
In particular, $U+V$   defines a radially symmetric solution of  the Liouville problem:
\ber\left\{\ba{lll}
  \ds -\Delta(U+V)=2|x|^{2(N_1+N_2)}\re^{U+V}\quad \text{in}\quad  \mathbb{R}^2,\\[2mm]
  \ds \int_{\mathbb{R}^2} |x|^{2(N_1+N_2)}\re^{U+V}\ud x <\infty,
 \ea\right.  \label{r26}
\eer
 which    satisfies: $\max(U+V)=U(0)+V(0)=0.$

Therefore, by the classification result of \cite{prta}, we get the following  explicit expression:
 \ber
  U(r)+V(r)=\ln\left(\frac{1}{\left(1+\frac{r^{2(N_1+N_2+1)}}{4(N_1+N_2+1)^2}\right)^2}\right),\label{r31}
 \eer
and
\ber
 \int_0^\infty r^{2(N_1+N_2)+1}\re^{U(r)+V(r)}\ud r =2(N_1+N_2+1).\label{r32}
\eer
 In addition,  by using (the radial expression of) \eqref{r25'} and \eqref{r31} we readily check that $U$ and $V$ admit
 the same logarithmic growth at $\infty$. As a consequence,  $U-V$  defines  a  bounded harmonic  function in $\mathbb{R}^2$,   which vanishes at the origin,
  in other words $U-V \equiv0$,   and we obtain:
\ber
 U(r)=\ln\left(\frac{1}{1+\frac{r^{2(N_1+N_2+1)}}{4(N_1+N_2+1)^2}}\right)=V(r),\quad \forall\, r\ge0. \label{r34}
\eer
In particular, we observe that:  $\frac{t_n}{s_n}\to 0$.

We can use those information to see that, for any $0<\vep<\frac14\min\{N_1, N_2\}$ sufficiently small, there exist $n_\vep\in \mathbb{N}$ and  $R_\vep\gg1$ such that, for $n\ge n_\vep$ and $r\ge R_\vep$,
we have:
\ber
\tilde{F}_{1,n}(r)&\equiv&\int_0^{r}t^{2(N_1+N_2)+1}\re^{\tilde{U}_n(t)+\tilde{V}_n(t)}\ud t+\vep_{2,n}\int_0^rt^{2N_2}\re^{\tilde{V}_n(t)}\ud t\nm\\
&\ge&\int_0^{r}t^{2(N_1+N_2)+1}\re^{\tilde{U}_n(t)+\tilde{V}_n(t)}\ud t> 2(N_1+N_2+1-\vep).\label{r33}
\eer
As a consequence, we can check that for  $n\ge n_\vep$, the function:
\ber
r^{2(N_1+1+\vep)}\re^{\tilde{U}_n(r)} \quad \text{is decreasing}\quad  \forall\, r\ge R_\vep,\label{r35}
\eer
indeed:
\ber
 \frac{\ud }{\ud r}\left(r^{2(N_1+1+\vep)}\re^{\tilde{U}_n(r)}\right)&=&r^{2(N_1+\vep)+1}\re^{\tilde{U}_n(r)}\left[2(N_1+1+\vep)-\tilde{F}_{1,n}(r)\right]\nm\\
  &\le&r^{2(N_1+\vep)+1}\re^{\tilde{U}_n(r)}\left[2(N_1+1+\vep)-2(N_1+N_2+1-\vep)\right]\nm\\
  &=&-2r^{2(N_1+\vep)+1}\re^{\tilde{U}_n(r)}(N_2-2\vep)<0,\quad \forall \, r\ge R_\vep.\label{r37}
\eer

Therefore, for $ r\ge R_\vep$ we find:
\ber
 r^{2(N_1+1+\vep)}\re^{\tilde{U}_n(r)}\le R_\vep^{2(N_1+1+\vep)}\re^{\tilde{U}_n(R_\vep)}\to \frac{R_\vep^{2(N_1+1+\vep)}}{1+\frac{R_\vep^{(N_1+N_2+1)}}{4(N_1+N_2+1)^2}},  \quad \text{as}\quad n\to \infty;\label{r38}
\eer
and so, for a suitably constant $C_\vep>0$,  we have:
\berr
r^{2(N_1+1+\vep)}\re^{\tilde{U}_n(r)}\le C_\vep,\quad \forall\, r\ge R_\vep,\,\, \forall \, n\ge n_\vep.
\eerr
 As a consequence,
\ber
\int_0^\infty r^{2N_1+1}\re^{\tilde{U}_n(r)}\ud r<C,  \label{r40}
\eer
for suitable $C>0$,  and we  conclude   the important fact,
\ber
\int_0^\infty r^{2N_1+1}\re^{U_n(r)}\ud r=\vep_{1, n}\int_0^\infty r^{2N_1+1}\re^{\tilde{U}_n(r)}\to 0, \,\,\text{as}\,\, n\to \infty.\label{r41}
\eer

 Similarly, we can check that (by taking $n_\vep$ larger if necessary)
\ber
r^{2(N_2+1+\vep)}\re^{\tilde{V}_n(r)} \quad \text{is decreasing}\quad  \forall r\ge R_\vep, \,\,\forall\,n\ge n_\vep, \label{r42}
\eer
 simply  by computing:
\ber
 &&\frac{\ud }{\ud r}(r^{2(N_2+1+\vep)}\re^{\tilde{V}_n(r)})\nm\\
 &&=r^{2(N_2+\vep)+1}\re^{\tilde{V}_n(r)}(2(N_2+1+\vep)-\tilde{F}_{2,n}(r))\nm\\
 &&\equiv r^{2(N_2+\vep)+1}\re^{\tilde{V}_n(r)}\left[2(N_2+1+\vep)-\int_0^rt^{2(N_1+N_2)}\re^{\tilde{U}_n(t)+\tilde{V}_n(t)}\ud r+\vep_{1,n}\int_0^rt^{2N_1}\re^{\tilde{U}_n(t)}\ud t\right]\nm\\
  &&\le r^{2(N_2+\vep)+1}\re^{\tilde{V}_n(r)}[2(N_2+1+\vep)- 2(N_1+N_2+1-\vep)+{\rm o}(1)]\nm\\
  &&\le -2r^{2(N_2+\vep)+1}\re^{\tilde{V}_n(r)}(N_1-2\vep+{\rm o}(1))<0,\quad \text{as}\quad n\to \infty. \label{r44}
\eer
 Therefore, exactly as above we can show that,
\ber
\int_0^\infty r^{2N_2+1}\re^{V_n(r)}\ud r=\vep_{2, n}\int_0^\infty r^{2N_2+1}\re^{\tilde{V}_n(r)}\ud r \to 0\,\,\text{as}\,\, n\to \infty.\label{r45}
\eer

 With this information, we can finally show that (by taking $n_\vep$ larger if necessary)
\ber
 r^{2(N_1+N_2+1+\vep)}\re^{\tilde{U}_n+\tilde{V}_n} \quad \text{is decreasing}, \,\, \forall\,r\ge R_\vep,\, \, \forall\, n\ge n_\vep, \label{r46}
\eer
  as  we have:
\ber
 &&\frac{\ud}{\ud r}\left[r^{2(N_1+N_2+1+\vep)}\re^{\tilde{U}_n(r)+\tilde{V}_n(r)}\right]\nm\\
 &&=r^{2(N_1+N_2+\vep)+1}\re^{\tilde{U}_n(r)+\tilde{V}_n(r)}\left[2(N_1+N_2+1+\vep)-\tilde{F}_{1,n}(r)-\tilde{F}_{2,n}(r)\right]\nm\\
 &&\le r^{2(N_1+N_2+\vep)+1}\re^{\tilde{U}_n(r)+\tilde{V}_n(r)}\bigg[2(N_1+N_2+1+\vep)-2\int_0^rt^{2(N_1+N_2)+1}\re^{\tilde{U}_n(t)+\tilde{V}_n(t)}+{\rm o}(1)\bigg] \nm\\
 &&\le-2r^{2(N_1+N_2+\vep)+1}\re^{\tilde{U}_n(r)+\tilde{V}_n(r)}(N_1+N_2-3\vep+{\rm o}(1))<0, \quad\text{as}\,\, n\to\infty.\label{r47}
\eer
So, we can use the convergence:  $(\tilde{U}_n+\tilde{V}_n)(R_\vep)\to (U+V)(R_\vep)$ as $n\to \infty$,  to  obtain as above  a suitable constant $C_\vep>0$
 such that,  for $n\ge n_\vep$ and  $r\ge R_\vep$,
\ber
  r^{2(N_1+N_2)+1}\re^{\tilde{U}_n+\tilde{V}_n} \le \frac{C_\vep}{r^{1+2\vep}},
\eer
 which allows us to  conclude that,  $\text{as}\, n\to\infty$:
\ber
  \int_0^\infty r^{2(N_1+N_2)+1}\re^{\tilde{U}_n(r)+\tilde{V}_n(r)}\ud r\to \int_0^\infty r^{2(N_1+N_2)+1}\re^{U(r)+V(r)}\ud r=2(N_1+N_2+1).\label{6.60}
\eer

Clearly,  \eqref{r41}, \eqref{r45}   and \eqref{6.60} imply that:
\ber
F_{2, n}(\infty)\to 2(N_1+N_2+1),\quad \text{as}\quad n\to \infty,
\eer
and  \eqref{5.1a} is established.  \hfill $\square$\\[2mm]

{\bf The proof of Theorem \ref{thfr}}.

In view of Theorem \ref{thbp}, for any fixed  constant $L\in \mathbb{R}$,   we take the initial data in the form:
\ber
 (\alpha_1, \alpha_2)\equiv\left(\alpha,\,  \frac{(2N_2+1)\alpha-L}{2N_1+1}\right), \quad \alpha\in \mathbb{R}, \label{id}
\eer
 and denote by $F_{2,\alpha}(\infty)$ the expression in  \eqref{s00}  corresponding to the unique solution $(U_\alpha, V_\alpha)$ of \eqref{s1} with initial data
 specified in \eqref{id}. To simplify notations, we do not emphasize the dependence  on $L$.

By virtue of Corollary \ref{corl62},   there exists $a_0$ such that, for any $\alpha<a_0$,
\ber
  F_{2,\alpha}(\infty)<2(N_2+1)\quad (\text{or equivalently}\,\,  F_{1,\alpha}(\infty)= \infty).
\eer

On the other hand,  by   Theorem \ref{thbp}, we  know  that,
\ber
   F_{2,\alpha}(\infty)\to 2(N_1+N_2+1),  \quad \text{as}\quad \alpha\to \infty.
\eer

Therefore, by the continuity of $F_{2,\alpha}(\infty)$ (or $F_{1,\alpha}(\infty)$) with respect to $\alpha$, for any $\gamma\in \big(2(N_2+1),\, 2(N_1+N_2+1)\big)$, there exists (at least) a parameter
 $\alpha_\gamma=\alpha_\gamma(L)$ such that the Cauchy problem \eqref{s1} with initial data given by \eqref{id} with $\alpha=\alpha_\gamma$
 admits a solution $(U_{\alpha_\gamma}, V_{\alpha_\gamma})$  satisfies the  integrability condition \eqref{z1} and  $F_{2,\alpha_\gamma}(\infty)=\gamma$. In other words,  we get  a radial  solution  of  \eqref{ss11}
 in the form $(u_{\alpha_\gamma}, v_{\alpha_\gamma})=(U_{\alpha_\gamma}+2N_1\ln r, V_{\alpha_\gamma}+2N_2\ln r)$.

Since for $L_1\neq L_2$ the  corresponding  set of initial data  satisfies:
 \[\left(\alpha_\gamma(L_1),\,  \frac{(2N_2+1)\alpha_\gamma(L_1)-L_1}{2N_1+1}\right)\neq \left(\alpha_\gamma(L_2),\,  \frac{(2N_2+1)\alpha_\gamma(L_2)-L_2}{2N_1+1}\right)\]
  in this way we have obtained a 1-parameter  family of different solutions of \eqref{ss11},  all satisfying:
 $F_2(\infty)=\gamma\in \big(2(N_2+1),\, 2(N_1+N_2+1)\big)$.

Furthermore, from the version \eqref{ss15} of the  Pohozaev identity,  it follows that:
as $\alpha\to \infty$,  if $F_{2,\alpha}(\infty)\to 2(N_1+N_2+1)$, then $\int_0^\infty r^{2N_1+1}\re^{U_\alpha}\ud r\to 0$,   $\int_0^\infty r^{2N_2+1}\re^{V_\alpha}\ud r\to 0$ and  $F_{1,\alpha}(\infty)\to 2(N_1+N_2+1)$.

As a consequence, if we fix $\gamma\in\big(2(N_1+N_2+1),\, \infty\big)$,  we obtain the same  conclusion as above  with suitable  $\alpha_\gamma\in \mathbb{R}$: $F_{1,\alpha_\gamma}(\infty)=\gamma$.
 Similarly,  if we set:
  \ber
  E_\alpha=N(N-1)\sigma^3k\int_0^\infty\left(r^{2N_1+1}\re^{U_\alpha}+r^{2N_2+1}\re^{V_\alpha}\right)\ud r,
  \eer
then we see that: $E_\alpha\to 0$ as $\alpha\to \infty$, while $E_\alpha=\infty$ for $\alpha<a_0$, and for every $E>0$ the above  conclusion follows
  with suitable $\alpha_E\in \mathbb{R}:$ $E_{\alpha_E}=E$,  and  the existence part of Theorem \ref{thfr} is established.

Finally,  by virtue of \eqref{s43a} and \eqref{z2},  we check that:  $\beta_1\equiv\frac12F_1(\infty)-N_1>N_2+1$ and $1<\beta_2\equiv\frac12F_2(\infty)-N_2<N_1+1$. Moreover,
  to obtain \eqref{s68}--\eqref{s70}, we recall that,
    \ber
   \lim\limits_{r\to \infty}ru'(r)=-2\beta_1, \quad \lim\limits_{r\to \infty}rv'(r)=-2\beta_2,
  \eer
  and so,  by taking the limit $r\to \infty$ in \eqref{e4'}--\eqref{e6'}, we have:
\ber
 &&2\int_0^\infty t\re^{u(t)+v(t)}\ud t=4(N_1+1)(N_2+1)-4(\beta_1-1)(\beta_2-1),\label{s72}\\
 &&2\int_0^\infty t\re^{u(t)}\ud t=4N_1(N_2+1)-4\beta_1(\beta_2-1),\label{s73}\\
 &&2\int_0^\infty t\re^{v(t)}\ud t=4\beta_2(\beta_1-1)-4N_2(N_1+1), \label{s74}
 \eer
from which \eqref{s68}--\eqref{s70} easily  follow and the proof of Theorem \ref{thfr} is complete. \hfill $\square$\\[2mm]

{\bf The proof of Theorem \ref{th1} and Theorem \ref{th2}}.

By virtue of the results established  above about problem \eqref{ss11}, we can easily  obtain the proof of Theorem \ref{th1} and Theorem \ref{th2}.

Indeed, by recalling \eqref{i48}--\eqref{i49} (considered with the lower sign)  and \eqref{re3'}, we see that,  via \eqref{re1}--\eqref{re2'},  the following
holds  for the fluxes:
 \ber
  &&\frac{1}{2\pi}\Phi_1=\int_0^\infty f^{(1)}_{12}(r)r\ud r=\frac{8\pi^2}{k^2}\int_0^\infty |\phi_2(r)|^2\left(|\phi_1(r)|^2+c^2\right)r\ud r\nm\\
  &&=2\sigma^2\int_0^\infty\re^{v(r)}\left(\re^{u(r)}+1\right)r\ud r=\frac12F_1(\infty), \label{s75}\\
  &&\frac{1}{2\pi}\Phi_2=\int_0^\infty f^{(2)}_{12}(r)r\ud r=\frac{8\pi^2}{k^2}\int_0^\infty |\phi_1(r)|^2\left(|\phi_2(r)|^2-c^2\right)r\ud r\nm\\
  &&=2\sigma^2\int_0^\infty\re^{u(r)}\left(\re^{v(r)}-1\right)r\ud r=\frac12F_2(\infty), \label{s76}
 \eer
 where  the last identity in \eqref{s75} and \eqref{s76} takes into account  the scaling \eqref{2.78a} and \eqref{i66}.
 Similarly, for the total energy we have:
 \berr
  E=N(N-1)\sigma^3k\int_0^\infty \left(\re^{u(r)}+\re^{v(r)}\right) r\ud r.
 \eerr

At this point, still by keeping in mind  \eqref{re1}--\eqref{re3},   we see that Theorem \ref{th1} and Theorem \ref{th2}  follow  directly by   \eqref{ss17} (or \eqref{s43a}),
Theorem \ref{thfr} and \eqref{de1}--\eqref{de2}.

\small{

}
\end{document}